\pdfoutput=1

\documentclass[onefignum,onetabnum]{siamart190516}
\usepackage{bm}
\usepackage{siunitx}
\usepackage[paper,cmrgreekup]{pdef}
\usepackage{algorithmic}
\usepackage{subfigure}
\usepackage{multirow}
\usepackage{makecell}
\usepackage{adjustbox}
\usepackage{booktabs}
\usepackage{pgf}
\usepackage{tikz}
\usetikzlibrary{calc,arrows.meta,decorations.pathmorphing,backgrounds}

\headers{Multi-Level Fine-Tuning}{Z. Li, Y. Fan, and L. Ying}

\title{
Multi-Level Fine-Tuning: Closing Generalization Gaps \\ In Approximation of Solution Maps \\ Under a Limited Budget for Training Data
\thanks{Submitted to the editors DATE.
\funding{Zhihan Li is partially supported by the elite undergraduate training program of the School of Mathematical Sciences at Peking University. Lexing Ying is partially supported by the National Science Foundation under award DMS-1818449 and by the U.S. Department of Energy, Office of Science, Office of Advanced Scientific Computing Research, Scientific Discovery through Advanced Computing (SciDAC) program.}}
}
\author{
Zhihan Li\thanks{School of Mathematical Sciences, Peking University, Beijing 100871, China (\email{zhli16@pku.edu.cn}).}
\and Yuwei Fan\thanks{Department of Mathematics, Stanford University, Stanford, CA 94305 (\email{ywfan1989@gmail.com}).}
\and Lexing Ying\thanks{Department of Mathematics and Institute for Computational and Mathematical Engineering, Stanford University, Stanford, CA 94305 (\email{lexing@stanford.edu}).}
}
\date{March 20, 2020}

\ifpdf
\hypersetup{
pdftitle={Multi-Level Fine-Tuning: Closing Generalization Gaps in Approximation of Solution Maps under a Limited Budget for Training Data},
pdfauthor={Z. Li, Y. Fan, and L. Ying}
}
\fi

\DeclareMathOperator\ope{\mathrm{E}}
\DeclareMathOperator\opeh{\hat{\mathrm{E}}}

\DeclareMathOperator\opsupp{\mathrm{supp}}

\begin{document}

\maketitle

\begin{abstract}
In scientific machine learning, regression networks have been recently applied to approximate solution maps (e.g., potential-ground state map of Schr\"odinger equation).
In this paper, we aim to reduce the generalization error without spending more time in generating training samples.
However, to reduce the generalization error, the regression network needs to be fit on a large number of training samples (e.g., a collection of potential-ground state pairs).
The training samples can be produced by running numerical solvers, which takes much time in many applications.
In this paper, we aim to reduce the generalization error without spending more time in generating training samples.
Inspired by few-shot learning techniques, we develop the Multi-Level Fine-Tuning algorithm by introducing levels of training: we first train the regression network on samples generated at the coarsest grid and then successively fine-tune the network on samples generated at finer grids.
Within the same amount of time, numerical solvers generate more samples on coarse grids than on fine grids.
We demonstrate a significant reduction of generalization error in numerical experiments on challenging problems with oscillations, discontinuities, or rough coefficients.
Further analysis can be conducted in the Neural Tangent Kernel regime and we provide practical estimators to the generalization error.
The number of training samples at different levels can be optimized for the smallest estimated generalization error under the constraint of budget for training data.
The optimized distribution of budget over levels provides practical guidance with theoretical insight as in the celebrated Multi-Level Monte Carlo algorithm.
\end{abstract}



\section{Introduction} \label{Sec:Intro}

\subsection{Background} \label{Sec:Back}

\subsubsection{Approximating solution maps with regression networks}

Contemporary machine learning techniques, especially deep neural networks, have been recently introduced to scientific computing tasks.
For problems involving partial differential equations (PDEs), neural networks can be utilized as universal approximators to the solution function in order to solve a PDE directly \cite{lagaris_artificial_1998,perdikaris_nonlinear_2017,e_deep_2018,khoo_solving_2018,raissi_deep_2018}.
In this setting, to solve the PDE $ \mathcal{N} u = 0 $ where $\mathcal{N}$ is a differential operator and $u$ is the solution function on domain $ \Omega \subseteq \mathbb{R}^d $, one searches for a neural network $\mathsf{NN}$ which inputs the coordinate $ \bm{x} \in \mathbb{R}^d $ and outputs an approximation to the solution function $ \mathsf{NN} \rbr{\bm{x}} \approx u \rbr{\bm{x}} $.
Another setting tackles the problem of parametric PDE $ \mathcal{N}_v u_v = 0 $ with variable parameter $v$ (e.g., coefficients, initial values, source terms).
In this setting, one may train a neural network to approximate the solution map $ v \mapsto u_v $ (the parameter-solution map) in order to apply to various parameters \cite{khoo_solving_2020}.
For example, we consider the one-dimensional non-linear Schr\"odinger equation on $ \Omega = \sbr{ 0, 1 } $ with a periodic boundary condition and a fixed dispersion coefficient $ \beta > 0 $ \cite{bagnato_bose--einstein_2015,fan_multiscale_2019}:
\begin{gather}
\label{Eq:NLSE1}
-\Delta u \rbr{x}  + v \rbr{x} u \rbr{x} + \beta u^3 \rbr{x} = E u \rbr{x}, \quad x \in \sbr{ 0, 1 }, \\
\label{Eq:NLSE2}
\int_{\sbr{ 0, 1 }} u^2 \rbr{x} \sd x = 1,  \quad \int_{\sbr{ 0, 1 }} u \rbr{x} \sd x > 0.
\end{gather}
We are interested in the solution map from potential $v$ (parameter) to ground state $u_v$ (the function $ u_v = u $ satisfying the equations with the smallest energy $E$, solution).
Hence, we train a neural network which inputs potentials and outputs ground states to approximate the potential-ground state map $ v \mapsto u_v $.
Besides, in graphics, one may approximate the action of Poisson solver, namely the solution map of $ -\Delta u_v = v $ from source term $v$ (parameter) to $u_v$ (solution) for faster Eulerian fluid simulation \cite{tompson_accelerating_2017,kim_deep_2019}.
For inverse problems, one may approximate the regularized inverse operator which maps observations (parameter) to reconstructions (solution) in a data-driven context \cite{khoo_switchnet_2019,arridge_solving_2019,fan_solving_2019,fan_solving_2020,gilton_neumann_2020}.
Similar settings of approximating solution maps can also be seen in signal processing \cite{sun_solving_2003}, molecular dynamics \cite{zhang_deep_2018}, model reduction \cite{lee_model_2020}, and operator compression \cite{fan_multiscale_2019-1,fan_multiscale_2019,fan_bcr-net_2019}.

In this paper, we focus on the setting of approximating solution maps.
For discretization, one typically chooses a fine grid of interest (in contrast to the coarse grid introduced later)
and discretizes the solution map as a fine-grid numerical solver.
Since the fine grid can be chosen up to practical considerations, main problem of approximating the solution map is transformed into approximating the fine-grid solver.
In order to harvest a neural network which approximates the fine-grid solver, one may fit a regression network (a neural network used as a regressor) on parameter-solution pairs generated by invoking the fine-grid solver.
For example, in the non-linear Schr\"odinger equation problem, we independently sample a collection of potentials $ \cbr{v_m}{}_{ m = 1 }^M $ from a distribution $\mathcal{D}$ as functions on grid with grid step $ h = 1 / 320 $.
We then compute the corresponding ground states $ u_m = u_{v_m} $ by running a gradient flow solver \cite{bao_computing_2004} on the grid.
With the parameter-solution pairs $ \bcbr{}{ \rbr{ v_m, u_m } }{}_{ m = 1 }^M $, we initialize and train a neural network $\mathsf{NN}$ to fit the training samples $ \rbr{ v_m, u_m } $.

\subsubsection{Tradeoff between generalization and training data generation}

The error of regression network in approximating the fine-grid solver consists of two parts: (1) the training error on $ \rbr{ v_m, u_m } $ and (2) the generalization error when applying the regression network on previously unseen samples (parameters) $ v \sim \mathcal{D} $.
In many cases, the neural network fits exactly at training samples because of over-parameterization \cite{belkin_understand_2018,liang_just_2020,allen-zhu_convergence_2019}, and hence the generalization error is of major concern.

The procedure described above to train a regression network lies in the framework of empirical risk minimization: we fit a statistical model (regression network) in a hypothesis space (e.g., functions representable by a neural network with bounded weights) by minimizing the empirical risk (training error) on training data, and then we evaluate the model on testing data in the hope of a low population risk (testing error).
Generalization gap (generalization error) is the difference between the population risk and the empirical risk \cite{mohri_foundations_2018}.
In particular, we have the following theorem bounding the generalization gap using Rademacher complexity \cite{mohri_foundations_2018}.
\begin{theorem} \label{Thm:RadeMar}
Let $\mathcal{L}$ be a family of functions from a set $\mathcal{V}$ to $ \sbr{ 0, B } $ and $ \cbr{v_m}{}_{ m = 1 }^M $ be $M$ independent random samples drawn from a distribution $\mathcal{D}$ on $\mathcal{V}$.
Then, for any positive $\delta$, with probability at least $ 1 - \delta $ on sampling $ \cbr{v_m}{}_{ m = 1 }^M $, it holds for every $ \ell \in \mathcal{L} $ that
\begin{equation}
\ope \ell \rbr{v} - \opeh \ell \rbr{v} \le 2 \hat{\mathfrak{R}}_v \rbr{\mathcal{L}} + 3 B \sqrt{\frac{ \log 1 / \delta }{ 2 M }}.
\end{equation}
\end{theorem}
Here $\mathcal{L}$ is the loss class (composition of statistical models and the loss function),
\begin{equation}
\ope \ell \rbr{v} = \ope_{ v \sim \mathcal{D} } \ell \rbr{v}, \quad \opeh \ell \rbr{v} = \frac{1}{M} \sum_{ m = 1 }^M \ell \rbr{v_m}
\end{equation}
are the population risk and the empirical risk respectively, and
\begin{equation} \label{Eq:Rade}
\hat{\mathfrak{R}}_v \rbr{\mathcal{L}} = \frac{1}{M} \ope_{\sigma} \sup_{ \ell \in \mathcal{L} }\sum_{ m = 1 }^M \sigma_m \ell \rbr{v_m}
\end{equation}
is the empirical Rademacher complexity, where $\sigma_m$ are independent Bernoulli-like random variables which take the value $ \pm 1 $ with equal probability.

Rademacher complexity delineates the complexity of hypothesis space and relates to the generalization gap.
In various cases, $ \hat{\mathfrak{R}}_v \rbr{\mathcal{L}} $ decays in speed $ O \frbr{ 1 / \sqrt{M} } $ (e.g., for linear class or kernel class \cite{bartlett_rademacher_2001,cortes_learning_2013}).
However, if $M$ is small, the Rademacher complexity can be poorly bounded, and so can the generalization gap.

As a result, in order to accurately approximate the fine-grid solver (and the solution map), a sufficient number of training samples are in need.
However, numerical solvers may suffer from heavy scaling laws in computational complexity with respect to grid step due to stability, convergence, or randomness issues.
Hence, generating a large number of training samples with the fine-grid solver may take a long time.
For example, in \cite{fan_solving_2019}, it takes about 50 seconds to solve a radiative transfer equation when generating a single parameter-solution pair.
As a result, the generation of the whole dataset, which consists of \num{10240} samples of such pairs, takes days and turns out to be much slower and much more expensive than the training of regression network.
We can observe the tradeoff between generalization and training data generation in this example.
Hence, it comes the question whether it is possible to reduce the generalization error in approximating the fine-grid solver without spending more time in generating training samples.
In other words, we aim to close the generalization gap under a limited budget for training data.

\subsubsection{Multi-Level Fine-Tuning}

Few-shot learning methods have been developed to mitigate the shortage of training samples.
Prior knowledge is required to either augment training samples, reduce the hypothesis space, or reach a better parameterization \cite{wang_generalizing_2020}.
One particularly popular method is fine-tuning, which was initially proposed for transfer learning tasks \cite{wang_generalizing_2020}.
In this approach, a model is first fit on some other tasks called source tasks with a large amount of data.
The model is then refined using a similar fitting procedure on a few-shot task called target task.
The intermediate model, as a good starting point for the target task, contains some prior knowledge about the source tasks, which may be transferred to help the few-shot target task.

For image recognition tasks, models pre-trained on large-scale datasets like ImageNet \cite{russakovsky_imagenet_2015} \cite{huh_what_2016} serve as excellent starting points for fine-tuning.
However, fine-grid numerical solvers to be approximated are problem-dependent.
It is difficult to find a pre-trained model well-suited for multiple different problems in scientific machine learning.
Meanwhile, according to multi-scale methods, we may capture rough shapes of solutions on a coarse grid and then refine details on the fine grid \cite{e_heterogenous_2003}.
The resemblance between fine-tuning techniques and multi-scale methods motivates our Multi-Level Fine-Tuning (MLFT) algorithm.
Since numerical solvers on coarser grids usually occupy less time due to the scaling laws, we may generate a great more of training samples with a coarse-grid numerical solver.
The source task of MLFT is to fit the regression network on training samples generated at the coarse grid, approximating the coarse-grid solver.
Because of the large number of coarse-grid samples, we may obtain a good approximation to the coarse-grid solver.
From the perspective of multi-scale methods, approximation to the coarse-grid solver provides macroscopic information about the fine-grid solver (and the solution map) which we want to approximate ultimately.
However, even if we have a good approximation to the coarse-grid solver, we still need to compensate for the difference between coarse-grid and fine-grid solvers.
Therefore, in the target task of MLFT, we fine-tune the regression network on a few training samples generated at the fine grid, approximating the fine-grid solver.
The fine-tuning step transfers information from the source task to the target task and hence reduces the generalization error in approximating the fine-grid solver.
In words of multi-scale methods, microscopic information gets refined in the fine-tuning step.

\begin{figure}[htbp]
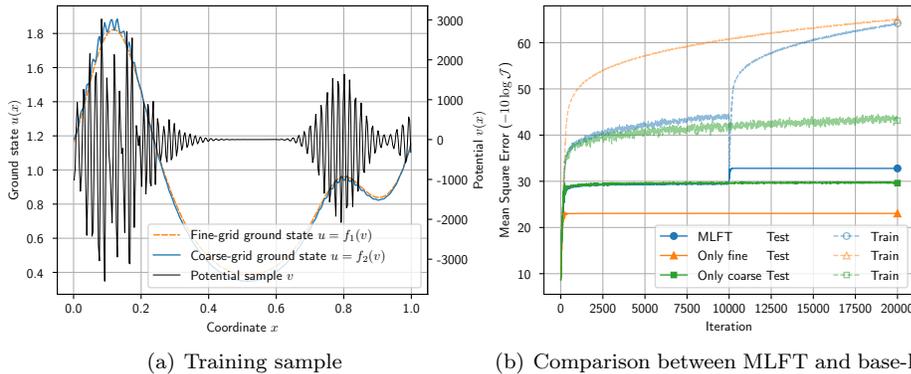

\centering
\subfigure[Training sample]{\draftbox{\trimbox{2pt 2pt -4pt 16pt}{\scalebox{0.5}{\input{Figures/NLSE1DFullDetFull.pgf}}}}}
\subfigure[Comparison between MLFT and base-lines]{\draftbox{\trimbox{4pt 2pt 4pt 16pt}{\scalebox{0.5}{\input{Figures/NLSE1DWideBigLoss.pgf}}}}}
\caption
{
The non-linear Schr\"odinger equation problem.
Since it is free to choose the finest grid of interest, we choose a fine grid with grid step $ h_2 = 1 / 320 $ and turn to approximate the fine-grid numerical solver $f_2$ with a regression network.
For MLFT, we introduce a coarse grid with grid step $ h_1 = 1 / 40 $ and a coarse-grid numerical solver $f_1$.
We plot a training sample generated by the two numerical solvers in (a).
In the comparison of loss curves (b), setting ``MLFT'' trains the regression network on $ M_1 = 1024 $ coarse-grid samples and then fine-tune on $ M_2 = 16 $ fine-grid samples, ``Only coarse'' trains on $ M_1 = 2048 $ coarse-grid samples, and ``Only fine'' trains on $ M_2 = 32 $ fine-grid samples.
The networks are tested on fine-grid samples.
We observe that the final testing error of MLFT is smaller than using either only coarse-grid or only fine-grid samples.
}
\label{Fig:Comp}
\end{figure}

In essence, the algorithm of two-level MLFT first trains the regression network on samples generated at the coarse grid, and then fine-tunes the the network on samples generate at the fine grid.
The weights and biases learned at the previous level of coarse grid are refined at the level of fine grid.
We take the problem of non-linear Schr\"odinger equation as an example.
As mentioned before, we aim to approximate the fine-grid solver with grid step $ h_2 = h = 1 / 320 $.
We assume that we generate $M_2$ training samples at the fine grid and each one takes time $t_2$.
We also introduce a coarse-grid solver with grid step $ h_1 = 1 / 40 $ to generate $M_1$ training samples at the coarse grid, each of which takes time $t_1$.
We observe from experiment that $ t_2 = 64 t_1 $.
The two-level MLFT algorithm consists of two stages: (1) we first initialize a regression network $\mathsf{NN}$ and fit it on the $M_1$ coarse-grid training samples and (2) we then fine-tune the regression network $\mathsf{NN}$ to fit the $M_2$ fine-grid training samples.
Given the budget of time $ T = 32 t_2 $ for generating training samples, we have the constraint $ M_1 t_1 + M_2 t_2 = T $ or $ M_1 / 64 + M_2 = 32 $.
We expect a reduction of generalization error using the two-level MLFT algorithm (e.g., $ M_1 = 1024 $ coarse-grid samples and $ M_2 = 16 $ fine-grid samples) compared to directly fitting only on $ M_2 = T / t_2 = 32 $ fine-grid samples.
We also expect the error of MLFT to be lower than the difference between fine-grid and coarse-grid solvers, or equivalently the testing error on fine-grid samples of a regression network trained only with $ M_1 = T / t_1 = 2048 $ coarse-grid samples.
The comparison is shown in \cref{Fig:Comp}.

By choosing a sequence of increasingly finer grids with grid step $ h_1 > h_2 > \cdots > h_L $, the two-level framework can be generalized to $L$ levels: we first train the regression network on samples generated at the coarsest grid ($h_1$) and then successively fine-tune the network on samples generated at finer grids ($h_2$, and then $h_3$, and so on till $h_L$).
We assume that $M_l$ training samples are generated by invoking the numerical solver on the grid with grid step $h_l$ and each one takes time $t_l$ on average, where $ t_1 < t_2 < \cdots < t_L $.
If an estimator to the generalization error is available in the form $ \hat{g}_L \rbr{ M_1, M_2, \cdots, M_L } $, we may optimize for the smallest estimated generalization error under the constraint of budget of time $T$ for generating training samples:
\begin{equation}
\label{Eq:Prog}
\begin{aligned}
\text{minimize} & & & \hat{g}_L \rbr{ M_1, M_2, \cdots, M_L }, \\
\text{subject to} & & & M_1 t_1 + M_2 t_2 + \cdots + M_L t_L = T, \\
\text{with respect to} & & & M_1, M_2, \cdots, M_L > 0.
\end{aligned}
\end{equation}
This budget distribution problem provides practical guidance to generate training samples at coarse and fine grids which invokes the coarse-grid and fine-grid numerical solvers respectively.
In other words, it gives insights on distributing the budget for training data over levels, as in the Multi-Level Monte Carlo (MLMC) algorithm \cite{giles_multilevel_2008}.

\subsection{Related work}

Structures of neural networks have inspired extensive study about hierarchical, multi-grid, and multi-level methods in machine learning, especially for convolutional and residual networks \cite{reed_parallel_2017,gao_learning_2018,fan_multiscale_2019-1,fan_multiscale_2019,juncai_mgnet_2019}.
In comparison, we focus on scales of grids on which numerical solvers are invoked and training samples are generated.
The main goal of MLFT is to reduce the generalization error in approximating the finest-grid solver (and the solution map).

Haber et al. considered an algebraic multi-grid method through the lens of optimal control \cite{haber_learning_2018}.
With designed restriction and interpolation procedures, the neural network can be transformed between scales of grids during training.
Varying the depths of networks can also be understood by temporally refining the optimal control problem.
In contrast, we have coarse-grid and fine-grid training samples, while only a single neural network is used with the architecture unchanged.
To match the size of training samples to the input and output shapes of the regression network, we apply spatial restriction and interpolation operators when generating training samples but not during training or fine-tuning the regression network.
Besides, the problems we are considering for scientific computing suffer from scaling laws of time
when generating training samples on different grids with numerical solvers.
On the contrary, classification and segmentation problems in machine learning need manual labeling, the cost of which generally do not depend on the resolution.

Very recently, Lye et al. also noticed the tradeoff between generalization and training data generation, and they developed the Multi-Level Machine Learning Monte Carlo (ML2MC) algorithm \cite{lye_multi-level_2020}.
Telescoping as in MLMC \cite{giles_multilevel_2008}, they train one regression network to approximate the coarsest-grid solver and multiple networks to approximate the difference between adjacent fine-grid and coarse-grid solvers.
The final model to approximate the finest-grid solver is the sum of all networks.
The generalization errors also add together, similar to the variance in MLMC.
In comparison, our initial motivation for MLFT follows from few-shot learning techniques.
We deploy only one neural network and fine-tune between collections of samples generated at different grids.
Since each level of fine-tuning reduces the distance between the current level solver and the intermediate network, generalization errors in previous levels get corrected in later levels.
This results in a form of generalization error different from ML2MC.
Additionally, \cite{lye_multi-level_2020} mainly considers the map from parameters to scalar observables.
In this paper, we tackle the problem of approximating numerical solvers (and the solution map), the input (parameter) and output (solution) of which are functions on grids and lie in high-dimensional spaces.

Another topic related to MLFT is multi-fidelity modeling \cite{perdikaris_multi-fidelity_2015,perdikaris_nonlinear_2017,peherstorfer_survey_2018}.
The setting of multi-fidelity models applies to ours: a combination of high-fidelity and low-fidelity models are accessible with different tradeoffs between efficiency and accuracy.
The high-fidelity and low-fidelity models correspond to fine-grid and coarse-grid solvers respectively.
In comparison, MLFT only accesses models (invokes solvers to generate training samples) of different levels subsequently and is more coarse-grained in terms of combining high-fidelity and low-fidelity models.
Meanwhile, we better leverage the nature of neural networks as parameterized statistical models which makes fine-tuning possible.
Recent progress in deep learning theory also provides us opportunities to find practical estimators to the generalization error.

\subsection{Contribution}

We summarize our contribution in this paper as follows.
\begin{partlist}
\item We identify the problem of reducing generalization error of regression networks in approximation of solution maps when the budget for generating training samples is limited (\cref{Sec:Intro}).
\item We design the Multi-Level Fine-Tuning (MLFT) algorithm to reduce generalization error with inspiration from few-shot learning techniques (\cref{Sec:AlgoSec}).
\item We perform analysis under the Neural Tangent Kernel (NTK) regime and construct practical estimators to the generalization error, which further provides guidance to distribute the budget for training data over levels (\cref{Sec:Ana}).
\item We show the reduction of generalization error with MLFT in experiments and demonstrate optimizing the number of training samples over levels (\cref{Sec:Num}).
\end{partlist}

\section{Algorithm} \label{Sec:AlgoSec}

We proceed to present our Multi-Level Fine-Tuning (MLFT) algorithm in this section with details.
As pointed out in \cref{Sec:Back}, the key idea of MLFT is to first train on samples generated at the coarsest grid and then successively fine-tune on samples generated at finer grids.
We introduce the definition of levels and the procedure to generate training samples of different levels in \cref{Sec:Level}.
We then explain the MLFT algorithm in \cref{Sec:Algo} together with the performance evaluation procedure of generalization error.
We compare MLFT with Multi-Level Machine Learning Monte Carlo (ML2MC) \cite{lye_multi-level_2020} in \cref{Sec:MLML}.

\subsection{Levels and data} \label{Sec:Level}

We tackle the problem of approximating solution maps in this paper.
For simplicity, we consider problems on a $d$-dimensional ($ d = 1, 2, 3 $) domain $ \Omega = \sbr{ 0, 1 }^d $ with a periodic boundary condition.
We consider the parametric PDE problem $ \mathcal{N}_v u_v = 0 $ where the variable parameter $v$ lies in a function space $ \mathcal{V} $.
We aim to approximate the non-linear solution map $ F : \mathcal{V} \to \mathcal{U} $ which maps parameter $ v \in \mathcal{V} $ to solution $ u_v \in \mathcal{U} $, where $ \mathcal{U} $ is the function space of solutions.
In the example of non-linear Schr\"odinger equation, the solution map $F$ is the potential-ground state map and we can set $ \mathcal{V} = \mathcal{U} = C \rbr{\Omega} $.

To numerically discretize the solution map $F$, we choose a finest grid of interest $ \Omega_L = \cbr{ i h_L : 0 \le i < N_L }^d $ of $N_L^d$ evenly spaced nodes with grid step $ h_L = 1 / N_L $.
Since we are considering regular grids, we identify functions on the grid as multi-dimensional arrays in $\mathbb{R}^{N_L^d}$.
By discretizing the PDE $ \mathcal{N}_v u_v = 0 $ on grid $\Omega_L$, we discretize the solution map $F$ to be the finest-grid solver $F_L$.
Since it is free to choose the finest grid step $h_L$ according to practical considerations, the core task is to approximate the finest-grid solver $F_L$ (discretized) instead of the solution map $F$ (continuous).
With parameters and solutions discretized as functions on grid $\Omega_L$, we assume the spaces of discretized parameters and solutions to be $ \mathcal{V}_L, \mathcal{U}_L \subseteq \mathbb{R}^{N_L^d} $ respectively.
In this way, the finest-grid numerical solver on grid $\Omega_L$ is $ F_L : \mathcal{V}_L \to \mathcal{U}_L $ which maps discretized parameters (functions on grid $\Omega_L$) to discretized solutions (functions on grid $\Omega_L$), as a discretization to the solution map $ F : \mathcal{V} \to \mathcal{U} $.
For example, in the one-dimensional non-linear Schr\"odinger equation problem, we consider the finest-grid solver on the finest grid with grid step $ h_2 = 1 / 320 $.
We may directly take the spaces of discretized potentials and ground states to be $ \mathcal{V}_2 = \mathcal{U}_2 = \mathbb{R}^{320} $, and the corresponding gradient flow ground state solver \cite{bao_computing_2004} is represented as $ F_2 : \mathcal{V}_2 \to \mathcal{U}_2 $.

As explained in \cref{Sec:Back}, the limited budget for generating training samples results in large generalization error.
To reduce the generalization error, MLFT introduces a series of coarser grids.
The regression network is fit on training samples generated at the coarser grids before finally fine-tuning on samples generated at the finest grid $\Omega_L$.
Formally, we choose a sequence of $L$ increasingly finer grid steps $ h_1 > h_2 > \cdots > h_{ L - 1 } > h_L $.
We apply similar numerical discretization to the parametric PDE $ \mathcal{N}_v u_v = 0 $ on the grid $ \Omega_l = \cbr{ i h_l : 0 \le i < N_l }^d $ with $ N_l = 1 / h_l $ for $ 1 \le l \le L - 1 $ as in the case of the finest grid $\Omega_L$.
We assume the spaces of potentials and ground states on grid $\Omega_l$ to be $ \mathcal{V}_l, \mathcal{U}_l \subseteq \mathbb{R}^{N_l^d} $ respectively, and the numerical solver working on grid $\Omega_l$ to be $ F_l : \mathcal{V}_l \to \mathcal{U}_l $.
For the example of non-linear Schr\"odinger equation, we use $ L = 2 $ levels and introduce a coarse grid with grid step $ h_1 = 1 / 40 $ for MLFT.
In this case, we set $ \mathcal{V}_1 = \mathcal{U}_1 = \mathbb{R}^{40} $ and the coarse-grid solver is represented as $ F_1 : \mathcal{V}_1 \to \mathcal{U}_1 $.

To approximate the finest-grid solver $F_L$, we introduce a regression network $\mathsf{NN}$.
By identifying multi-dimenisonal arrays as functions on grids, the network input and output functions on some grid.
In order to capture details on the finest grid of interest $\Omega_L$, we design the network to input and output functions on the finest grid $\Omega_L$, as denoted by $ \mathsf{NN} : \mathcal{V}_L \to \mathcal{U}_L $.
In other words, the network $\mathsf{NN}$ works on grid $\Omega_L$.
In the example of non-linear Schr\"odinger equation, we adopt a MNN-$\mathcal{H}$ network \cite{fan_multiscale_2019} to approximate the finest-grid solver $ F_2 : \mathcal{V}_2 = \mathbb{R}^{320} \to \mathcal{U}_2 = \mathbb{R}^{320} $.
The input and output dimensions of the network are both $ \sbr{ \text{\texttt{batch\_size}}, 320 } $, since the finest grid of interest $\Omega_2$ has $ N_2 = 320 $ nodes.

However, for $ 1 \le l \le L - 1 $, the coarse-grid solver $F_l$ works on the coarse grid $\Omega_l$ and is not compatible with the input and output dimensions of the network.
To generate coarse-grid training samples as functions on grid $\Omega_L$, we introduce restriction operators $ R_{ L \to l } : \mathcal{V}_L \to \mathcal{V}_l $ and interpolation operators $ I_{ l \to L } : \mathcal{U}_l \to \mathcal{U}_L $.
We transform the coarse-grid solver $ F_l : \mathcal{V}_l \to \mathcal{U}_l $ to
\begin{equation}
f_l = I_{ l \to L } \circ F_l \circ R_{ L \to l } : \mathcal{V}_L \to \mathcal{U}_L
\end{equation}
which works on the finest grid $\Omega_L$.
We train the network with collections of training samples which are all functions on the finest grid $\Omega_L$ but are generated by different level-$l$ solvers $f_l$.
For the non-linear Schr\"odinger equation problem, the coarse-grid solver $ F_1 : \mathcal{V}_1 = \mathbb{R}^{40} \to \mathcal{U}_1 = \mathbb{R}^{40} $ generates potential-ground state pairs with spatial resolution $ N_1 = 40 $, incompatible with the spatial resolution of network $ N_2 = 320 $.
Hence, we introduce a Fourier restriction operator $ R_{ 2 \to 1 } : \mathcal{V}_2 = \mathbb{R}^{320} \to \mathcal{V}_1 = \mathbb{R}^{40} $ and a bicubic interpolation operator $ I_{ 1 \to 2 } : \mathcal{U}_1 = \mathbb{R}^{40} \to \mathcal{U}_2 = \mathbb{R}^{320} $.
The transformed function of coarse-grid solver is $ f_1 = I_{ 1 \to 2 } \circ F_1 \circ R_{ 2 \to 1 } : \mathcal{V}_2 = \mathbb{R}^{320} \to \mathcal{U}_2 = \mathbb{R}^{320} $.
We use $f_1$ instead of $F_1$ to generate coarse-grid samples for training.
For notation convenience, we denote $ f_L = F_L $, $ R_{ L \to L } = \mathrm{id}_{\mathcal{V}_L} $, and $ I_{ L \to L } = \mathrm{id}_{\mathcal{U}_L} $.

\begin{figure}[htbp]
\centering
\draftbox{\trimbox{6pt 6pt 6pt 6pt}{\input{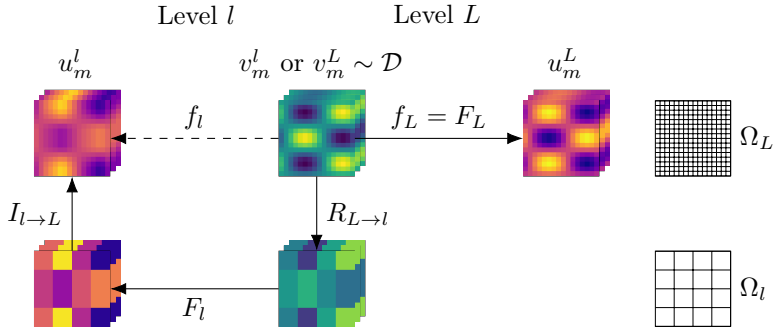}}}
\caption
{
Schematic illustration on generating samples at levels $l$ and $L$ where $ 1 \le l \le L - 1 $.
For level $l$, we draw parameter samples from $\mathcal{D}$ on grid $\Omega_L$, restrict the parameters to grid $\Omega_l$, invokes the coarse-grid solver $F_l$ for solutions on grid $\Omega_l$, and interpolate the solutions to grid $\Omega_L$.
}
\label{Fig:GenImg}
\end{figure}

We describe the procedure to generate samples at different levels.
According to the parameters of interest, we choose a probability distribution $\mathcal{D}$ on $\mathcal{V}_L$ to sample discretized parameters.
At level $l$ where $ 1 \le l \le L $, we independently draw $M_l$ parameter samples $ \bcbr{}{v^l_m}{}_{ m = 1 }^{M_l} $ from $\mathcal{D}$ on grid $\Omega_L$.
We then restrict the parameters to grid $\Omega_l$, invoke the coarse-grid solver $F_l$ for solutions on grid $\Omega_l$, and interpolate the solutions to grid $\Omega_L$.
Equivalently, we compute $ u^l_m = f_l \brbr{}{v^l_m} = I_{ l \to L } \circ F_l \circ R_{ L \to l } \brbr{}{v^l_m} $, as depicted in \cref{Fig:GenImg}.
For finding the optimized number of training samples over levels, we denote the average time to generate a sample at level $l$ by evaluating $f_l$ to be $t_l$.

\subsection{Multi-Level Fine-Tuning algorithm} \label{Sec:Algo}

\begin{figure}[htbp]
\centering
\draftbox{\trimbox{6pt 6pt 6pt 6pt}{\begin{tikzpicture}
\input{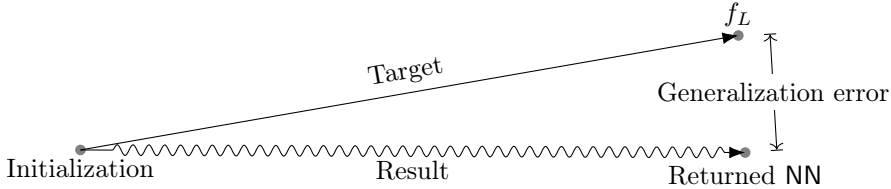}
\node [below] at (f0) {Initialization};
\node [above] at (f4) {$f_L$};
\node [below] at (fh4) {\begin{tabular}{c} Returned 
$\mathsf{NN}$ \end{tabular}};
\fill [gray] (f0) circle (2pt);
\fill [gray] (f4) circle (2pt);
\fill [gray] (fh4) circle (2pt);
\draw [->,arrow] (f0) -- node [sloped,above] {Target} (f4);
\draw [->,arrow,decorate,wave] (f0) -- node [sloped,below] {Result} (fh4);
\draw [|<->|] ($(f4)!12pt!90:(fh4)$) -- node [fill=white] {Generalization error} ($(fh4)!12pt!270:(f4)$);
\tikzbox
\end{tikzpicture}%
}}
\caption
{
Schematic illustration of \cref{Alg:TrainSing} (single-level training at level $ l = L $) in the function space.
In this illustration, points represent functions in $ \mathcal{V}_L \to \mathcal{U}_L $: we aim to approximate the finest-grid solver $f_L$ (a function from $\mathcal{V}_L$ to $\mathcal{U}_L$) by a regression network $\mathsf{NN}$ (a parameterized function from $\mathcal{V}_L$ to $\mathcal{U}_L$).
We target to compensate for the difference between the initialized network and the finest-grid solver $f_L$ (``Target''), and during training the network moves in the function space (``Result'').
When the training error $g^{\mathrm{train}}_L$ vanishes, difference between the trained regression network $\mathsf{NN}$ and the finest-grid solver $f_L$ is the generalization error $ g_L = g^{\mathrm{test}}_L $.
}
\label{Fig:ScheSing}
\end{figure}

We describe our Multi-Level Fine-Tuning (MLFT) algorithm.
In order to approximate the solution map $ F : \mathcal{V} \to \mathcal{U} $, we deploy a regression neural network $ \mathsf{NN} : \mathcal{V}_L \to \mathcal{U}_L $ to approximate the finest-grid numerical solver $ f_L : \mathcal{V}_L \to \mathcal{U}_L $.
With only one level $ l = L $, we fit the neural network on training samples generated at level $l$ as described in \cref{Sec:Level}.
We use the Mean Square Error (MSE), or equivalently the squared $L^2$ norm on grid $\Omega_L$ for this training (fitting) process as in \cref{Alg:TrainSing}.
We introduce the schematic illustration in the function space $ \mathcal{V}_L \to \mathcal{U}_L $ in \cref{Fig:ScheSing}.

\begin{algorithm}[htbp]
\caption{Train a regression network at a single level $l$}
\label{Alg:TrainSing}
\begin{algorithmic}
\STATE{Generate $M_l$ training samples $ \bcbr{}{ \brbr{}{ v^l_m, u^l_m } }{}_{ m = 1 }^{M_l} $ with $ u^l_m = f_l \brbr{}{v^l_m} $ at level $l$}
\STATE{Initialize a neural network $\mathsf{NN}$ which works on grid $\Omega_L$}
\STATE{Fit $\mathsf{NN}$ on $ \bcbr{}{ \brbr{}{ v^l_m, u^l_m } }{}_{ m = 1 }^{M_l} $}
by minimizing the Mean Square Error (MSE)
\begin{equation} \label{Eq:MSE}
\mathcal{J}_l = \frac{1}{M_l} \sum_{ m = 1 }^{M_l} \norm{ u^l_m - \mathsf{NN} \brbr{}{v^l_m} }_{L^2}^2
\end{equation}
\RETURN Trained regression neural network $\mathsf{NN}$
\end{algorithmic}
\end{algorithm}

The single-level training of a regression network at level $ l = L $ is the most direct approach to approximate the finest-grid solver $f_L$.
Single-level training at level $ l = L $ serves as a base-line for comparison.
As mentioned in \cref{Sec:Back}, since neural networks usually fit exactly at the training samples, the generalization error in approximating the finest-grid solver $f_L$ is of major concern.
Hence, we evaluate the generalization error $g_L$ in approximating the finest-grid solver $f_L$ for performance comparison:
\begin{equation} \label{Eq:GLDef}
g_L = g_L^{\mathrm{test}} - g_L^{\mathrm{train}}, \quad \\
\begin{cases}
g_L^{\mathrm{test}} = \ope_{ v \sim \mathcal{D} } \norm{ f_L \rbr{v} - \mathsf{NN} \rbr{v} }_2, \\
g_L^{\mathrm{train}} = \sum_{ m = 1 }^{M_L} \norm{ u^L_m - \mathsf{NN} \brbr{}{v^L_m} }_2 / M_L.
\end{cases}
\end{equation}
Here we use the vector 2-norm (instead of $L^2$-norm on grid $\Omega_L$) in consistence with the notation in \cref{Sec:Ana}.
We expect $ g_L^{\mathrm{train}} \ll g_L^{\mathrm{test}} $ and $ g_L \approx g_L^{\mathrm{test}} $.
Another base-line for comparison is the single-level training only with samples generated at a coarse grid, for example, the single-level training at level $ l = L - 1 $.
In this case, we evaluate the testing error for performance comparison since we do not train at level $L$:
\begin{equation} \label{Eq:GLNorm}
g_L = g_L^{\mathrm{test}} = \ope_{ v \sim \mathcal{D} }\norm{ f_L \rbr{v} - \mathsf{NN} \rbr{v} }_2.
\end{equation}
The regression network approximates the coarse-grid solver $ f_{ L - 1 } $ instead of the finest-grid solver $f_L$, so we expect $ g_L \gtrsim e_L $ where $e_L$ is the difference in between:
\begin{equation}
e_L = \ope_{ v \sim \mathcal{D} } \norm{ f_L \rbr{v} - f_{ L - 1 } \rbr{v} }_2.
\end{equation}

\begin{figure}[htbp]
\centering
\draftbox{\trimbox{6pt 6pt 6pt 6pt}{\begin{tikzpicture}
\input{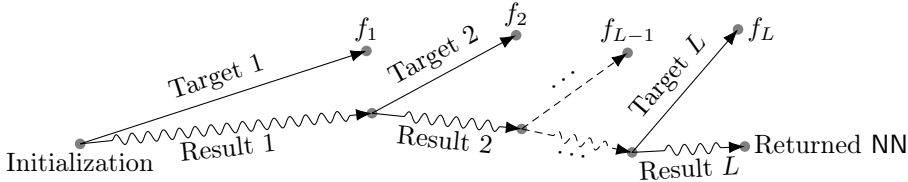}
\node [below] at (f0) {Initialization};
\node [above] at (f1) {$f_1$};
\node [above] at (f2) {$f_2$};
\node [above] at (f3) {$ f_{ L - 1 } $};
\node [right] at (f4) {$f_L$};
\node [right] at (fh4) {\hspace{-6pt}\begin{tabular}{c} Returned 
$\mathsf{NN}$ \end{tabular}};
\fill [gray] (f0) circle (2pt);
\fill [gray] (f1) circle (2pt);
\fill [gray] (f2) circle (2pt);
\fill [gray] (f3) circle (2pt);
\fill [gray] (f4) circle (2pt);
\fill [gray] (fh1) circle (2pt);
\fill [gray] (fh2) circle (2pt);
\fill [gray] (fh3) circle (2pt);
\fill [gray] (fh4) circle (2pt);
\draw [->,arrow] (f0) -- node [sloped,above] {Target $1$} (f1);
\draw [->,arrow] (fh1) -- node [sloped,above] {Target $2$} (f2);
\draw [->,arrow,densely dashed] (fh2) -- node [sloped,above] {$\cdots$} (f3);
\draw [->,arrow] (fh3) -- node [sloped,above] {Target $L$} (f4);
\draw [->,arrow,decorate,wave] (f0) -- node [sloped,below] {Result $1$} (fh1);
\draw [->,arrow,decorate,wave] (fh1) -- node [sloped,below] {Result $2$} (fh2);
\draw [->,arrow,decorate,wave,densely dashed] (fh2) -- node [sloped,below] {$\cdots$} (fh3);
\draw [->,arrow,decorate,wave] (fh3) -- node [sloped,below] {Result $L$} (fh4);
\tikzbox
\end{tikzpicture}%
}}
\caption
{
Schematic illustration of \cref{Alg:TrainMLFT} (MLFT) in the function space.
We mark the training targets (``Target $l$'') and the obtained results (``Result $l$'') at level $l$, as in \cref{Fig:ScheSing}.
}
\label{Fig:ScheMLFT}
\end{figure}

As mentioned in \cref{Sec:Back}, our MLFT algorithm first trains the regression neural network on samples generated at the coarsest grid and then successively fine-tunes the network on samples generated at finer grids.
The MLFT algorithm is described in \cref{Alg:TrainMLFT} together with the function space illustration in \cref{Fig:ScheMLFT}.
Similar to the single-level training, MSE \cref{Eq:MSE} is used as the loss function for the regression network to fit the training samples.
Since we deploy only one neural network $\mathsf{NN}$, there is only one initialization step.
In the following fine-tuning steps, we do not freeze parameters (weights and biases) of the neural network, nor modify the network architecture.
We have not observed apparent over-fitting phenomenon in numerical experiments approximating solution maps and numerical solvers (e.g., \cref{Fig:Comp,Fig:BurgErr,Fig:Ellip2}), so we do not apply extra regularization either.
The optimizer is not restart but we keep the momentum vector (for Momentum \cite{qian_momentum_1999}) or the estimation of moments (for Adam \cite{kingma_adam_2015}).
For simplicity, parameters of the optimizer like step sizes are retained during the whole training and fine-tuning process.

\begin{algorithm}[htbp]
\caption{Multi-Level Fine-Tune (MLFT) a regression network}
\label{Alg:TrainMLFT}
\begin{algorithmic}
\STATE{Generate $M_1$ samples $ \bcbr{}{ \brbr{}{ v^1_m, u^1_m } }{}_{ m = 1 }^{M_1} $ with $ u^1_m = f_1 \brbr{}{v^1_m} $}
\STATE{Initialize a neural network $\mathsf{NN}$ which works on grid $\Omega_L$}
\STATE{Fit $\mathsf{NN}$ on $ \bcbr{}{ \brbr{}{ v^1_m, u^1_m } }{}_{ m = 1 }^{M_1} $ by minimizing MSE}
\STATE{\hfill}\COMMENT{Train $\mathsf{NN}$ from initialization to $f_1$}
\FOR{$ l = 2 $ \TO $L$}
\STATE{Generate $M_l$ samples $ \bcbr{}{ \brbr{}{ v^l_m, u^l_m } }{}_{ m = 1 }^{M_l} $ with $ u^l_m = f_l \brbr{}{v^l_m} $}
\STATE{Fit $\mathsf{NN}$ on $ \bcbr{}{ \brbr{}{ v^l_m, u^l_m } }{}_{ m = 1 }^{M_l} $ by minimizing MSE}
\STATE{\hfill}\COMMENT{Fine-tune $\mathsf{NN}$ to $f_l$}
\ENDFOR
\RETURN Trained and fine-tuned regression neural network $\mathsf{NN}$
\end{algorithmic}
\end{algorithm}

\subsection{Multi-Level Machine Learning Monte Carlo algorithm} \label{Sec:MLML}

In \cite{lye_multi-level_2020}, Lye et al. proposed a multi-level algorithm for regression networks named Multi-Level Machine Learning Monte Carlo (ML2MC) with inspiration from Multi-Level Monte Carlo (MLMC) \cite{giles_multilevel_2008}.
In the paper \cite{lye_multi-level_2020}, the problem of approximating parameter-observable maps of parametric PDEs is considered, where the observables are scalars.
We make slight modifications to match our setting of approximating solution maps (parameter-solution maps), whose inputs and outputs are functions on grids.

\begin{figure}[htbp]
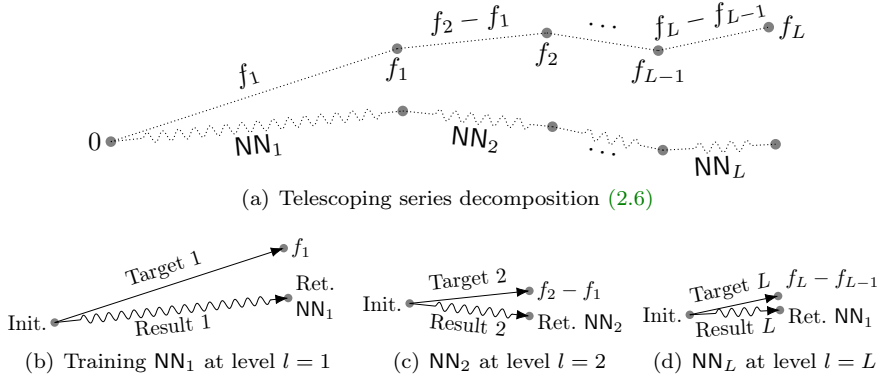

\centering
\subfigure[Telescoping series decomposition \cref{Eq:ML2MC}]{\draftbox{\trimbox{6pt 6pt 6pt 6pt}{\begin{tikzpicture}
\input{Schematics/Header.tikz}
\node [left] at (f0) {$0$};
\node [below] at (f1) {$f_1$};
\node [below] at (f2) {$f_2$};
\node [below] at (f3) {$ f_{ L - 1 } $};
\node [right] at (f4) {$f_L$};
\fill [gray] (f0) circle (2pt);
\fill [gray] (f1) circle (2pt);
\fill [gray] (f2) circle (2pt);
\fill [gray] (f3) circle (2pt);
\fill [gray] (f4) circle (2pt);
\fill [gray] (fh1) circle (2pt);
\fill [gray] (fh2) circle (2pt);
\fill [gray] (fh3) circle (2pt);
\fill [gray] (fh4) circle (2pt);
\draw [densely dotted] (f0) -- node [sloped,above] {$f_1$} (f1);
\draw [densely dotted] (f1) -- node [sloped,above] {$ f_2 - f_1 $} (f2);
\draw [densely dotted] (f2) -- node [sloped,above] {$\cdots$} (f3);
\draw [densely dotted] (f3) -- node [sloped,above] {$ f_L - f_{ L - 1 } $} (f4);
\draw [decorate,wave,densely dotted] (f0) -- node [sloped,below] {$\mathsf{NN}_1$} (fh1);
\draw [decorate,wave,densely dotted] (fh1) -- node [sloped,below] {$\mathsf{NN}_2$} (fh2);
\draw [decorate,wave,densely dotted] (fh2) -- node [sloped,below] {$\cdots$} (fh3);
\draw [decorate,wave,densely dotted] (fh3) -- node [sloped,below] {$\mathsf{NN}_L$} (fh4);
\tikzbox
\end{tikzpicture}%
}}}
\subfigure[Training $\mathsf{NN}_1$ at level $ l = 1 $]{\draftbox{\trimbox{0pt 3pt 0pt 3pt}{\scalebox{0.8}{\begin{tikzpicture}
\input{Schematics/Header.tikz}
\coordinate [label=left:Init.] (a1) at (0.0, 0.0);
\coordinate [label=right:$f_1$] (t1) at ($ (a1) + (f1) - (f0) $);
\coordinate [label={[label distance=-6pt]right:%
\begin{tabular}{c} Ret. \\
$\mathsf{NN}_1$ \end{tabular}}
] (r1) at ($ (a1) + (fh1) - (f0) $);
\fill [gray] (a1) circle (2pt);
\fill [gray] (t1) circle (2pt);
\fill [gray] (r1) circle (2pt);
\draw [->,arrow] (a1) -- node [sloped,above] {Target $1$} (t1);
\draw [->,arrow,decorate,wave] (a1) -- node [sloped,below] {Result $1$} (r1);
\tikzbox
\end{tikzpicture}%
}}}}
\subfigure[$\mathsf{NN}_2$ at level $ l = 2 $]{\draftbox{\trimbox{3pt 3pt 3pt 3pt}{\scalebox{0.8}{\begin{tikzpicture}
\input{Schematics/Header.tikz}
\coordinate [label=left:Init.] (a2) at (0.0, 0.0);
\coordinate [label=right:$ f_2 - f_1 $] (t2) at ($ (a2) + (f2) - (f1) $);
\coordinate [label={[label distance=-8pt]below right:
\begin{tabular}{c} Ret.\ 
$\mathsf{NN}_2$ \end{tabular}}
] (r2) at ($ (a2) + (fh2) - (fh1) $);
\fill [gray] (a2) circle (2pt);
\fill [gray] (t2) circle (2pt);
\fill [gray] (r2) circle (2pt);
\draw [->,arrow] (a2) -- node [sloped,above] {Target $2$} (t2);
\draw [->,arrow,decorate,wave] (a2) -- node [sloped,below] {Result $2$} (r2);
\tikzbox
\end{tikzpicture}%
}}}}
\subfigure[$\mathsf{NN}_L$ at level $ l = L $]{\draftbox{\trimbox{3pt 3pt 3pt 3pt}{\scalebox{0.8}{\begin{tikzpicture}
\input{Schematics/Header.tikz}
\coordinate [label=left:Init.] (a4) at (0.0, 0.0);
\coordinate [label=above right:$ f_L - f_{ L - 1 } $] (t4) at ($ (a4) + (f4) - (f3) $);
\coordinate [label={[label distance=-8pt]below right:
\begin{tabular}{c} Ret.\ 
$\mathsf{NN}_1$ \end{tabular}}
] (r4) at ($ (a4) + (fh4) - (fh3) $);
\fill [gray] (a4) circle (2pt);
\fill [gray] (t4) circle (2pt);
\fill [gray] (r4) circle (2pt);
\draw [->,arrow] (a4) -- node [sloped,above] {Target $L$} (t4);
\draw [->,arrow,decorate,wave] (a4) -- node [sloped,below] {Result $L$} (r4);
\tikzbox
\end{tikzpicture}%
}}}}
\caption
{
Schematic illustrations of \cref{Alg:TrainML2MC} (ML2MC) in the function space.
Function of the finest-grid solver $f_L$ is decomposed by the telescoping series \eqref{Eq:ML2MC} as shown in (a), and then $L$ neural networks $\mathsf{NN}_l$ are trained separately at level $l$ for $ 1 \le l \le L$, as shown in (b), (c), and (d).
}
\label{Fig:ScheML2MC}
\end{figure}

The algorithm of ML2MC constructs a telescoping series as in MLMC
\begin{equation} \label{Eq:ML2MC}
f_L = f_1 + \rbr{ f_2 - f_1 } + \rbr{ f_3 - f_2 } + \cdots + \rbr{ f_L - f_{ L - 1 } }.
\end{equation}
One regression network is used to approximate $f_1$ and other $ L - 1 $ networks to approximate $ f_l - f_{ l - 1 } $ for $ 2 \le l \le L $ respectively.
For an approximation to the finest-grid solver $f_L$, one sums the $L$ separately-trained networks together, as explained in \cref{Alg:TrainML2MC} and \cref{Fig:ScheML2MC}.

\begin{algorithm}[htbp]
\caption{Multi-Level Machine Learning Monte Carlo (ML2MC) for regression networks \cite{lye_multi-level_2020}}
\label{Alg:TrainML2MC}
\begin{algorithmic}
\STATE{Generate $M_1$ samples $ \bcbr{}{ \brbr{}{ v^1_m, u^1_m } }{}_{ m = 1 }^{M_1} $ with $ u^1_m = f_1 \brbr{}{v^1_m} $}
\STATE{Initialize a neural network $\mathsf{NN}_1$ which works on grid $\Omega_L$}
\STATE{Fit $\mathsf{NN}_1$ on $ \bcbr{}{ \brbr{}{ v^1_m, u^1_m } }{}_{ m = 1 }^{M_1} $ by minimizing MSE}
\STATE{\hfill}\COMMENT{Fit $\mathsf{NN}_1$ from initialization to $f_1$}
\FOR{$ l = 2 $ \TO $L$}
\STATE{Generate $M_l$ samples $ \bcbr{}{ \brbr{}{ v^l_m, u^l_m } }{}_{ m = 1 }^{M_l} $ with $ u^l_m = f_l \brbr{}{v^l_m} - f_{ l - 1 } \brbr{}{v^l_m} $}
\STATE{Initialize a neural network $\mathsf{NN}_l$ which works on grid $\Omega_L$}
\STATE{Fit $\mathsf{NN}_l$ on $ \bcbr{}{ \brbr{}{ v^l_m, u^l_m } }{}_{ m = 1 }^{M_l} $ by minimizing MSE}
\STATE{\hfill}\COMMENT{Fit $\mathsf{NN}_l$ from initialization to $ f_l - f_{ l - 1 } $}
\ENDFOR
\RETURN Sum of trained regression neural networks $ \sum_{ l = 1 }^L \mathsf{NN}_l $
\end{algorithmic}
\end{algorithm}

The main difference between MLFT (\cref{Alg:TrainMLFT}) and ML2MC (\cref{Alg:TrainML2MC}) is two-fold.
The first difference is that we fine-tune one regression network instead of summing several separately-trained networks.
We make use of the difference between parameterized model fitting and Monte Carlo estimation:
fine-tuning may correct generalization errors accumulated in previous levels, while variances of multiple independent estimators (and generalization errors of several separately-trained networks) add together.
As a result, the form of generalization error of MLFT is different from that of ML2MC, as will be analyzed in \cref{Sec:Ana}.
Another difference lies in the usage of neural networks.
In MLFT, we are able to confine all the computation into a single neural network thanks to the fine-tuning technique.
This saves graphic memory and avoids restarting the optimizer: modern first-order optimizers memorize historical information for acceleration \cite{qian_momentum_1999,kingma_adam_2015}.
Moreover, insights from curriculum learning \cite{weinshall_curriculum_2018} tell that we may gain much more.
If a parameterized statistical model is trained on easy tasks before moving on to more difficult ones, the training process on difficult tasks can get boosted significantly.
In our setting, approximating the finest-grid solver $f_L$ which contains more details can be more difficult than approximating coarse-grid solvers $f_l$ where $ 1 \le l \le L - 1 $.
In correspondence, we observe faster convergence in the fine-tuning steps of MLFT compared to ML2MC, as shown in experiments of \cref{Sec:Num,Fig:BurgErr}.

\section{Analysis} \label{Sec:Ana}

We analyze our algorithm of MLFT under assumptions in this section.
In the MLFT algorithm, at level $l$, we fit the regression network on training samples $ \bcbr{}{ \brbr{}{ v^l_m, u^l_m } }_{ m = 1 }^{M_l} $ where $ u^l_m = f_l \brbr{}{v^l_m} $ to approximate the solver $f_l$ where $ 1 \le l \le L $.
As a result, generalization errors are presented at the training or fine-tuning process at each level.
We consider bounds of the generalization error at each level citing results of Neural Tangent Kernel (NTK) and Rademacher complexity of kernel classes in \cref{Sec:NTK,Sec:GenNTK}.
We then chain the generalization errors at each level under guidance of empirical observations and construct \emph{a priori} error estimator $\hat{g}_L$ to the generalization error $g_L$ of MLFT in approximating the finest-grid solver in \cref{Sec:APri}.
We try to get rid of pessimistic estimations and extend to finite-width cases not solidly covered by the infinite-width NTK theory by introducing practical \emph{a posteriori} error estimator $\hat{g}_L$ by fitting coefficients into the form of generalization error in \cref{Sec:APost}.

\begin{figure}[htbp]
\centering
\draftbox{\trimbox{6pt 6pt 6pt 6pt}{\begin{tikzpicture}
\input{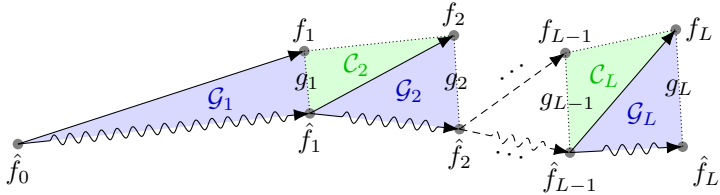}
\node [below] at (f0) {$\hat{f}_0$};
\node [above] at (f1) {$f_1$};
\node [above] at (f2) {$f_2$};
\node [above] at (f3) {$ f_{ L - 1 } $};
\node [right] at (f4) {$f_L$};
\node [below] at (fh1) {$\hat{f}_1$};
\node [below] at (fh2) {$\hat{f}_2$};
\node [below] at (fh3) {$ \hat{f}_{ L - 1 } $};
\node [below right] at (fh4) {$\hat{f}_L$};
\fill [fill=blue!15] (f0) -- (fh1) -- (f1) -- cycle;
\fill [fill=green!15] (fh1) -- (f2) -- (f1) -- cycle;
\fill [fill=blue!15] (fh1) -- (fh2) -- (f2) -- cycle;
\fill [fill=green!15] (fh3) -- (f4) -- (f3) -- cycle;
\fill [fill=blue!15] (fh3) -- (fh4) -- (f4) -- cycle;
\fill [gray] (f0) circle (2pt);
\fill [gray] (f1) circle (2pt);
\fill [gray] (f2) circle (2pt);
\fill [gray] (f3) circle (2pt);
\fill [gray] (f4) circle (2pt);
\fill [gray] (fh1) circle (2pt);
\fill [gray] (fh2) circle (2pt);
\fill [gray] (fh3) circle (2pt);
\fill [gray] (fh4) circle (2pt);
\draw [->,arrow] (f0) -- node [sloped,above] {} (f1);
\draw [->,arrow] (fh1) -- (f2);
\draw [->,arrow,densely dashed] (fh2) -- node [sloped,above,pos=0.6] {$\cdots$} (f3);
\draw [->,arrow] (fh3) -- (f4);
\draw [->,arrow,decorate,wave] (f0) -- (fh1);
\draw [->,arrow,decorate,wave] (fh1) -- (fh2);
\draw [->,arrow,decorate,wave,densely dashed] (fh2) -- node [sloped,below] {$\cdots$} (fh3);
\draw [->,arrow,decorate,wave] (fh3) -- (fh4);
\draw [densely dotted] (fh1) -- node 
{$g_1$} (f1);
\draw [densely dotted] (fh2) -- node 
{$g_2$} (f2);
\draw [densely dotted] (fh3) -- node 
{$ g_{ L - 1 } $} (f3);
\draw [densely dotted] (fh4) -- node 
{$g_L$} (f4);
\draw [densely dotted] (f1) -- node [sloped,above] {} (f2);
\draw [densely dotted] (f3) -- node [sloped,above,pos=0.6] {} (f4);
\node [text=blue!75!black] at (2.7, 0.6) {$\mathcal{G}_1$};
\node [text=green!75!black] at (4.5, 1.05) {$\mathcal{C}_2$};
\node [text=blue!75!black] at (5.2, 0.7) {$\mathcal{G}_2$};
\node [text=green!75!black] at (7.8, 0.9) {$\mathcal{C}_L$};
\node [text=blue!75!black] at (8.3, 0.4) {$\mathcal{G}_L$};
\tikzbox
\end{tikzpicture}%
}}
\caption
{
Schematic illustration of notations and the construction of error estimators (cf. \cref{Fig:ScheMLFT}).
For triangles {\color{blue!75!black}$\mathcal{G}_l$}, we estimate the Generalization error $g_l$ (i.e., $ f_l - \hat{f}_l $) from the initial difference at training samples $ u^l_m - \hat{f}_{ l - 1 } \brbr{}{v^l_m} $ (i.e.,  $ f_l - \hat{f}_{ l - 1 } $).
For triangles {\color{green!75!black}$\mathcal{C}_l$}, we Combine the generalization error $ g_{ l - 1 } $ (i.e., $ f_{ l - 1 } - \hat{f}_{ l - 1 } $) and the difference of training samples between levels $ u^l_m - f_{ l - 1 } \brbr{}{v^l_m} $ (i.e., $ f_l - f_{ l - 1 } $) to estimate the initial difference at training samples $ u^l_m - \hat{f}_{ l - 1 } \brbr{}{v^l_m} $ (i.e., $ f_l - \hat{f}_{ l - 1 } $).
By chaining the triangles, we obtain an error estimator to the generalization error $g_L$ in approximating the finest-grid solver $f_L$.
}
\label{Fig:NoteEst}
\end{figure}

In this section, we consider the training process in the function space $ \mathcal{V}_L \to \mathcal{U}_L $ as in \cite{jacot_neural_2018} and \cref{Fig:ScheSing,Fig:ScheMLFT,Fig:ScheML2MC}.
We denote the function of $\mathsf{NN}$ at initialization as $\hat{f}_0$, together with $\hat{f}_l$ for the intermediate model right after training or fine-tuning at level $l$.
In terms of functions, the level $l$ involves training or fine-tuning the regression network which is initially $ \hat{f}_{ l - 1 } $ to fit the target $f_l$ on training samples (``Target $l$'' in \cref{Fig:ScheMLFT}), and the network eventually moves from $ \hat{f}_{ l - 1 } $ to $\hat{f}_l$ (``Result $l$'' in \cref{Fig:ScheMLFT}).
We denote the generalization error of training or fine-tuning at level $l$ to be
\begin{equation} \label{Eq:GllDef}
g_l = g_l^{\mathrm{test}} - g_l^{\mathrm{train}}  = \ope_{ v \sim \mathcal{D} } \bnorm{\big}{ f_l \rbr{v} - \hat{f}_l \rbr{v} }_2 - \frac{1}{M_l} \sum_{ m = 1 }^{M_l} \bnorm{\big}{ u^l_m - \hat{f}_l \brbr{}{v^l_m} }_2.
\end{equation}
The key observation in introducing the function space illustrations is that in the Neural Tangent Kernel (NTK) regime, fitting the regression neural network on training samples converges to kernel ``ridgeless'' regression \cite{liang_just_2020}.
The kernel ``ridgeless'' regression is linear with respect to the dependent variable (usually referred to as $y$ in contrast to independent variables $x$).
In our case, it turns out that at level $l$ the increment $ \hat{f}_l - \hat{f}_{ l - 1 } $ (``Result $l$'' in \cref{Fig:ScheMLFT}) and the error $ f_l - \hat{f}_l $ in the function space only depend on the initial difference $ f_l - \hat{f}_{ l - 1 } $ (``Target $l$'' in \cref{Fig:ScheMLFT}) but not the starting point $ \hat{f}_{ l - 1 } $.
This is critical to our construction of estimators to the generalization error.
We briefly introduce the idea in \cref{Fig:NoteEst}.

\subsection{Neural Tangent Kernel} \label{Sec:NTK}

To understand the optimization and generalization of neural networks, one particular approach is to consider the infinite-width limit of neural networks in the function space \cite{jacot_neural_2018,lee_wide_2019}.
Under the Neural Tangent Kernel (NTK) parameterization and random unit Gaussian initialization, channels of a neural network behave like independent samples.
Hence, the law of large numbers can be cited, proving the convergence of the neural network to a Gaussian process.

Formally speaking, we denote the output of the network to be $ \mathsf{NN} \rbr{v} = f \rbr{ v; \theta } \in \mathbb{R}^{N_L^d} $ with input $v$ and parameter (weights and biases) $\theta$.
We consider NTK parameterization for $\mathsf{NN}$, which scales the output of each layer by $ 1 / \sqrt{C} $ where $C$ is the number of channels in the layer \cite{jacot_neural_2018}.
By recognizing the output $ f \rbr{ v; \theta } $ as a column vector, the Conjugate Kernel (CK, also known as Neural Network Gaussian Process or NNGP) is defined as
\begin{equation}
\Sigma \rbr{ v, v' } = \ope_{\theta} f \rbr{ v; \theta } f^{\mathsf{T}} \rbr{ v'; \theta } \in \mathbb{R}^{ N_L^d \times N_L^d },
\end{equation}
where the expectation is taken with respect to random unit Gaussian initialization of $\theta$.
When the numbers of channels $C$ go to infinity, $\Sigma$ converges and the function $ v \mapsto f \rbr{ v; \theta } $ (random because of random initialization of $\theta$) turns out to be the centered matrix-valued Gaussian process with covariance kernel $\Sigma$ \cite{jacot_neural_2018}.

Similarly, the Neural Tangent Kernel (NTK) is defined as
\begin{equation} \label{Eq:NTKDef}
\Theta \rbr{ v, v' } = \sum_{\theta} \frac{ \sd f \rbr{ v; \theta } }{ \sd \theta } \frac{ \sd f^{\mathsf{T}} \rbr{ v'; \theta } }{ \sd \theta } \in \mathbb{R}^{ N_L^d \times N_L^d },
\end{equation}
where the sum is taken over all entries of the parameter $\theta$.
As pointed out by \cite{jacot_neural_2018}, $\Theta$ converges almost surely at random unit Gaussian initialization of $\theta$ in the infinite-width limit.
In the case of gradient descent, NTK stays asymptotically constant during training \cite{jacot_neural_2018}.
A closed-form formula of the converged CK $\Sigma$ and NTK $\Theta$ is available for various neural network architectures (including dense layers, convolutional layers, and ReLU activation layers) \cite{arora_exact_2019}.
Software package has been developed to compute CK $\Sigma$ and NTK $\Theta$ both in the infinite-width limit and for a finite-width network \cite{novak_neural_2020}.

We summarize our assumptions for the following analysis.
We assume NTK parameterization and random unit Gaussian initialization of the neural network.
Moreover, we assume the network is infinitely wide so that we consider the NTK regime.
The influence of finite width is discussed in \cref{Supp:Sec:NTKConv}.
We assume that $\Theta$ is a symmetric positive definite kernel, which can be proved under assumptions (e.g., when norms of inputs are bounded both above and below, i.e., $ \Theta \rbr{1} $) \cite{jacot_neural_2018,cao_generalization_2020}.
In this case, $\Theta$ itself induces a Reproducing Kernel Hilbert Space (RKHS), which we denote by $\mathcal{H}$.
According to \cite{cao_generalization_2019}, because of the architecture of neural networks, we have $ \Sigma \prec \Theta $ and hence $ \hat{f}_0 \in \mathcal{H} $.
Additionally, we assume that the target at each level, namely the function of numerical solver satisfy $ f_l \in \mathcal{H} $.
This means that the functions $f_l$ are learnable under kernel ``ridgeless'' regression with NTK $\Theta$.
The learnability of certain functions with NTK has been be proved \cite{arora_fine-grained_2019}.

At level $l$, we train a regression network to fit $M_l$ training samples $ \bcbr{}{\brbr{}{ v^l_m, u^l_m }}_{ m = 1 }^M $ by minimizing MSE (equivalent to \cref{Eq:MSE} up to a constant)
\begin{equation} \label{Eq:MSEJ}
\mathcal{J}_l \sbr{f} = \frac{1}{M_l} \sum_{ m = 1 }^{M_l} \norm{ u^l_m - f \brbr{}{ v^l_m; \theta } }_2^2.
\end{equation}
According to the framework of NTK \cite{jacot_neural_2018}, the gradient descent dynamics on \cref{Eq:MSEJ} of an infinitely wide network in the function space $\mathcal{H}$ turns out to be a linear ordinary differential equation
\begin{equation} \label{Eq:LinODE}
\frac{ \sd f \rbr{v} }{ \sd t } = \sum_{\theta} \frac{ \sd f \rbr{ v; \theta } }{ \sd \theta } \frac{ \sd \theta }{ \sd t } = \frac{1}{M_l} \sum_{ m = 1 }^{M_l} \Theta \brbr{}{ v, v^l_m } \brbr{}{ u^l_m - f \brbr{}{v^l_m} }.
\end{equation}

Since the training or fine-tuning process at level $l$ fits the regression network to training samples generated by the level-$l$ solver $f_l$, we have $ u^l_m = f_l \brbr{}{v^l_m} $ and $ \nvbr{f}_{ t = 0 } = \hat{f}_{ l - 1 } $.
By introducing the Gram operator $ \Pi \in \mathcal{L} \rbr{\mathcal{H}} $, \cref{Eq:LinODE} turns out to be
\begin{equation} \label{Eq:GramODE}
\frac{ \sd f }{ \sd t } = \Pi \rbr{ f_l - f }, \quad \Pi \sbr{f} \rbr{v} = \frac{1}{M_l} \sum_{ m = 1 }^{M_l} \Theta \brbr{}{ v, v^l_m } f \brbr{}{v^l_m}.
\end{equation}
Since the Gram operator $\Pi$ is self-adjoint, positive semi-definite, and finite-rank in $\mathcal{H}$, solution $f$ of \eqref{Eq:GramODE} converges when $ t \rightarrow +\infty $.
We assume that we train the regression network for sufficiently long time at level $l$.
In this case, the infinite-time limit $ \nvbr{f}_{ t \rightarrow +\infty } = \hat{f}_l $.
Moreover, due to the assumption that NTK $\Theta$ is symmetric positive definite, the neural network indeed fits exactly at the training samples $ \bcbr{}{ \brbr{}{ v^l_m, u^l_m } }_{ m = 1 }^{M_l} $.
This can be summarized as the following theorem \cite{jacot_neural_2018}.
Specifically, the increment $ \hat{f}_l - \hat{f}_{ l - 1 } $ at level $l$ (``Result $l$'' in \cref{Fig:ScheMLFT}) only depends on the initial difference $ f_l - \hat{f}_{ l - 1 } $ (``Target $l$'' in \cref{Fig:ScheMLFT}) but not the specific function $ \hat{f}_{ l - 1 } $.
Thanks to this observation, we are able to construct error estimators to the generalization error as in \cref{Fig:NoteEst,Sec:APri}.

\begin{theorem} \label{Thm:Proj}
Assume the NTK $\Theta$ to be symmetric positive definite and $ \hat{f}_0, f_l \in \mathcal{H} $ for $ 1 \le l \le L $.
Let $\mathcal{M}_l$ be the space spanned by column spaces of $ \Theta \rbr{ \cdot, x^l_m } $ for $ 1 \le m \le M_l $ and $\mathrm{P}_{\mathcal{M}_l}$ be the orthogonal projection operator onto $\mathcal{M}_l$ in $\mathcal{H}$.
Then, for $ 1 \le l \le L $, it holds that
\begin{equation} \label{Eq:Proj}
\hat{f}_l - \hat{f}_{ l - 1 } = \mathrm{P}_{\mathcal{M}_l} \rbr{ f_l - \hat{f}_{ l - 1 } }.
\end{equation}
Specifically, $ \hat{f}_l \in \mathcal{H} $ for $ 1 \le l \le L $, and $ \hat{f}_l \brbr{}{v^l_m} = f_l \brbr{}{v^l_m} = u^l_m $ for $ 1 \le m \le M_l $.
\end{theorem}

\subsection{Generalization in the Neural Tangent Kernel regime} \label{Sec:GenNTK}

We then consider the generalization error of regression network in the NTK regime (infinite-width limit).
Generalization errors of neural networks have been extensively studied from the perspective of kernel methods \cite{belkin_understand_2018,arora_fine-grained_2019,cao_generalization_2019}.
Owing to the introduction of NTK, fitting a neural network is equivalent to performing kernel ``ridgeless'' regression with NTK $\Theta$.
As a result, the generalization error $g_l$ in approximating the level-$l$ solver function $f_l$ is related to the Rademacher complexity of kernel classes with NTK $\Theta$.

\subsubsection{Kernel ``ridgeless'' regression}

We have studied the increment $ \hat{f}_l - \hat{f}_{ l - 1 } $ at level $l$ in the function space in \cref{Thm:Proj}.
However, what we are really interested is the RKHS norm $ \bnorm{\big}{ \hat{f}_l - \hat{f}_{ l - 1 } }_{\mathcal{H}} $, which controls the size of hypothesis space and determines the Rademacher complexity.
The RKHS norm can be derived from the fact that fitting a neural network is equivalent to performing kernel ``ridgeless'' regression with NTK $\Theta$ \cite{jacot_neural_2018}.
(For a detailed exposition on kernel methods, see \cite{hofmann_kernel_2008,mohri_foundations_2018}.)

The kernel ``ridgeless'' regression is the limit of kernel ridge regression with diminishing ridge regularization.
The kernel ``ridgeless'' regression has been recently found to generalize \cite{belkin_understand_2018,liang_just_2020}.
In our setting, the result function $ \hat{f}_l = \hat{f}_{ l - 1 } + \mathrm{P}_{\mathcal{M}_l} \brbr{}{ f_l - \hat{f}_{ l - 1 } } $ is exactly the solution $\hat{f}^{\star}$ to the optimization problem of kernel ``ridgeless'' regression:
\begin{equation} \label{Eq:Sat}
\begin{aligned}
\text{minimize} & & & \bnorm{\big}{ \hat{f} - \hat{f}_{ l - 1 } }_{\mathcal{H}}, \\
\text{subject to} & & & \hat{f} \brbr{}{v^l_m} = f_l \brbr{}{v^l_m} = u^l_m, \quad 1 \le m \le M_l, \\
\text{with respect to} & & & \hat{f} \in \mathcal{H}.
\end{aligned}
\end{equation}

Denote the $ \mathbf{G}_l \in \mathbb{R}^{ M_l N_L^d \times M_l N_L^d } $ to be the Gram matrix and $ \bm{f}_{l'} \brbr{}{\bm{v}^l} \in \mathbb{R}^{ M_l N_L^d }$ (and also $ \hat{\bm{f}}_{l'} \brbr{}{\bm{v}^l} $) to be vectorized functions on training samples $v^l_m$:
\begin{equation} \label{Eq:Stack}
\mathbf{G}_l := \msbr{ \Theta \brbr{}{ v^l_1, v^l_1 } & \Theta \brbr{}{ v^l_1, v^l_2 } & \cdots & \Theta \brbr{}{ v^l_1, v^l_{M_l} } \\ \Theta \brbr{}{ v^l_2, v^l_1 } & \Theta \brbr{}{ v^l_2, v^l_2 } & \cdots & \Theta \brbr{}{ v^l_2, v^l_{M_l} } \\ \vdots & \vdots & \ddots & \vdots \\ \Theta \brbr{}{ v^l_{M_l}, v^l_1 } & \Theta \brbr{}{ v^l_{M_l}, v^l_2 } & \cdots & \Theta \brbr{}{ v^l_{M_l}, v^l_{M_l} } }, \quad
\bm{f}_{l'} \brbr{}{\bm{v}^l} := \msbr{ f_{l'} \brbr{}{v^l_1} \\ f_{l'} \brbr{}{v^l_2} \\ \vdots \\ f_{l'} \brbr{}{v^l_{M_l}} }.
\end{equation}
Using the kernel trick, the solution $ \hat{f}_l = \hat{f}^{\star} $ to \cref{Eq:Sat} can be written as (cf. \cref{Eq:Proj}) \cite{hofmann_kernel_2008}
\begin{equation}
\hat{f}_l = \hat{f}_{ l - 1 } + \sum_{ m = 1 }^{M_l} \Theta \brbr{}{ \cdot, v^l_m } \alpha_m, \quad \bm{\alpha} = \mathbf{G}_l^{-1} \brbr{\big}{ \bm{f}_l \brbr{}{\bm{v}^l} - \hat{\bm{f}}_{ l - 1 } \brbr{}{\bm{v}^l} },
\end{equation}
where $ \alpha_m \in \mathbb{R}^{N_L^d} $ are column vectors and $\bm{\alpha}$ is their vectorization as in \eqref{Eq:Stack}.
Abbreviating $ \pbr{ \bm{u}, \bm{w} }_{\mathbf{G}_l^{-1}} = \bm{w}^{\mathsf{T}} \mathbf{G}_l^{-1} \bm{u} $ and $ \norm{\bm{u}}_{\mathbf{G}_l^{-1}} = \sqrt{\vphantom{\pbr{}}\smash[b]{\pbr{ \bm{u}, \bm{u} }_{\mathbf{G}_l^{-1}} }} $, the RKHS norm is
\begin{equation}
\bnorm{\big}{ \hat{f}_l - \hat{f}_{ l - 1 } }_{\mathcal{H}} = \norm{\bm{\alpha}}_{\mathbf{G}_l} = \bnorm{\big}{ \bm{f}_l \brbr{}{\bm{v}^l} - \hat{\bm{f}}_{ l - 1 } \brbr{}{\bm{v}^l} }_{\mathbf{G}_l^{-1}}.
\end{equation}
This result is summarized as the following theorem \cite{hofmann_kernel_2008}.
\begin{theorem} \label{Thm:Norm}
Assume the NTK $\Theta$ to be symmetric positive definite and $ \hat{f}_0, f_l \in \mathcal{H} $ for $ 1 \le l \le L $.
Then, for $ 1 \le l \le L $,
\begin{equation}
\bnorm{\big}{ \hat{f}_l - \hat{f}_{ l - 1 } }_{\mathcal{H}} = \bnorm{\big}{ \bm{f}_l \brbr{}{\bm{v}^l} - \hat{\bm{f}}_{ l - 1 } \brbr{}{\bm{v}^l} }_{\mathbf{G}_l^{-1}}.
\end{equation}
\end{theorem}

\subsubsection{Rademacher complexity of kernel classes}

The complexity of hypothesis space of kernel classes can be delineated by Rademacher complexity \cite{mohri_foundations_2018}.
With the RKHS norm $ \bnorm{\big}{ \hat{f}_l - \hat{f}_{ l - 1 } }_{\mathcal{H}} $ given in \cref{Thm:Norm}, we are going to consider the Rademacher complexity of $\mathcal{H}$-balls centered at $ \hat{f}_{ l - 1 } \in \mathcal{H} $.
\begin{equation}
\mathcal{B} \brbr{\big}{ \hat{f}_{ l - 1 }, D } = \bcbr{\big}{ f \in \mathcal{H} : \bnorm{\big}{ f - \hat{f}_{ l - 1 } }_{\mathcal{H}} \le D }.
\end{equation}
We follow \cite{maurer_rademacher_2006,sindhwani_scalable_2013} to define the Rademacher complexity of vector-valued functions (in $\mathcal{H}$) and matrix-valued kernels (NTK $\Theta$).
Given samples $ \bcbr{}{v^l_m}{}_{ m = 1 }^{M_l} \subseteq \mathcal{V}_L $ and random testing vectors $ \bcbr{}{\nu_m}{}_{ m = 1 }^{M_l} $ in $\mathbb{R}^{N_L^d}$, we define the empirical Rademacher complexity of a vector-valued function class $\mathcal{F}$ from $\mathcal{V}_L$ to $\mathbb{R}^{N_L^d}$ to be
\begin{equation}
\hat{\mathfrak{R}}_{ v^l, \nu } \rbr{\mathcal{F}} := \frac{1}{M_l} \ope_{ \sigma, \nu } \sup_{ f \in \mathcal{F} } \sum_{ m = 1 }^{M_l} \sigma_m \nu_m^{\mathsf{T}} f \brbr{}{v^l_m}
\end{equation}
Here $\sigma_m$ are independent Bernoulli-like random variables which take the value $ \pm 1 $ with equal probability.
As a variant of the Rademacher complexity of scalar-valued kernel classes \cite{bartlett_rademacher_2001,mohri_foundations_2018}, we have the following theorem.
Detailed proof is given in \cref{Supp:Sec:Pf}.

\begin{theorem} \label{Thm:RadeKer}
Assume $ \bnorm{\big}{ \Theta \brbr{}{ v^l_m, v^l_m } }_2 \le R^2 $ for $ 1 \le m \le M_l $.
If the normalization condition $ \norm{\nu_m}_2 \equiv 1 $ and the independence condition $ \nu_m \perp \sigma_m $ are satisfied for $ 1 \le m \le M_l$, then it holds that
\begin{equation}
\hat{\mathfrak{R}}_{ v^l, \nu } \brbr{\big}{ \mathcal{B} \brbr{\big}{ \hat{f}_{ l - 1 }, D } } \le \frac{ D R }{\sqrt{M}}.
\end{equation}
\end{theorem}

We then consider the Rademacher complexity of loss classes, for adaptation to \cref{Thm:RadeMar}.
Since we aim to approximate the level-$l$ solver function $f_l$ as ground-truth at level $l$, we define the loss class of $ \mathcal{B} \brbr{\big}{ \hat{f}_{ l - 1 }, D } $ to be
\begin{equation}
\mathcal{L} \brbr{\big}{ \hat{f}_{ l - 1 }, D; f_l, B_l } = \bcbr{\big}{ v \mapsto \bnorm{\big}{ f_l \rbr{v} - f \rbr{v} }_2 \wedge B_l : f \in \mathcal{B} \brbr{\big}{ \hat{f}_{ l - 1 }, D } }.
\end{equation}
Here $B_l$ is a cutoff constant to bound the loss function.
We imitate \cite{mohri_foundations_2018} to derive the Rademacher complexity of loss classes of vector-valued functions in the following theorem.
The proof can be found in \cref{Supp:Sec:Pf}.

\begin{theorem} \label{Thm:RadeLoss}
Assume $ \bnorm{\big}{ \Theta \brbr{\big}{ v^l_m, v^l_m } }_2 \le R^2 $ for $ 1 \le m \le M_l $.
Then, it holds that
\begin{equation}
\hat{\mathfrak{R}}_{v^l} \brbr{\big}{ \mathcal{L} \brbr{\big}{ \hat{f}_{ l - 1 }, D; f_l, B_l } } \le \frac{ D R }{\sqrt{M}}.
\end{equation}
\end{theorem}

\subsubsection{Generalization of kernel classes}

We finally consider bounds of the generalization error in approximating the level-$l$ solver function using NTK kernel classes.
As an application of \cref{Thm:RadeMar,Thm:RadeLoss} and a variant of \cite{arora_fine-grained_2019}, we have the following theorem which is proved in \cref{Supp:Sec:Pf}.
Since the network fits exactly at the training samples as in \cref{Thm:Proj},
the term of training error is absent in \cref{Eq:Mar}.

\begin{theorem} \label{Thm:Mar}
Assume $ \norm{ \Theta \rbr{ v, v } }_2 \le R^2 $ for all $v$ in the support of $\mathcal{D}$.
Then, with probability $ 1 - \delta $ on generating training samples $ \bcbr{}{ \brbr{}{ v^l_m, u^l_m } }{}_{ m = 1 }^{M_l} $, it holds that
\begin{multline} \label{Eq:Mar}
\ope_{ v \sim \mathcal{D} } \brbr{\big}{ \bnorm{\big}{ f_l \rbr{v} - \hat{f}_l \rbr{v} }_2 \wedge B_l }
\\
\le \frac{ 2 R \bnorm{\big}{ \bm{f}_l \brbr{}{\bm{v}^l} - \hat{\bm{f}}_{ l - 1 } \brbr{}{\bm{v}^l} }_{\mathbf{G}_l^{-1}} }{\sqrt{M_l}} + \frac{ 2 B_l R }{\sqrt{M_l}} + 3 B_l \sqrt{\frac{ \log K_l / \delta }{ 2 M_l }},
\end{multline}
where $ K_l = \bgbr{\big}{ \bnorm{\big}{ f_l - \hat{f}_{ l - 1 } }_{\mathcal{H}} / B_l } $.
\end{theorem}

\subsection{\emph{A priori} error estimator} \label{Sec:APri}

We aim to find estimators to the generalization error $g_L$ in approximating the finest-grid solver $f_L$.
Since under assumptions the network fits exactly at the training samples at level $l$ (cf. \cref{Thm:Proj}), the generalization error \cref{Eq:GllDef} turns out to be $ g_l = \ope_{ v \sim \mathcal{D} } \bnorm{\big}{ f_l \rbr{v} - \hat{f}_l \rbr{v} }_2 $.
We have already mentioned the equivalence between neural network and kernel ``ridgeless'' regression with NTK $\Theta$ and considered the generalization error of NTK kernel classes in \cref{Thm:Mar}.
The remaining task is to apply \cref{Thm:Mar} in practical cases and chaining the generalization error as introduced in \cref{Fig:NoteEst}.
In this subsection, we aim to find the \emph{a priori} error estimator which can be computed before performing MLFT.
We provide empirical observations to justify our construction of the estimator.

\subsubsection{Estimating \texorpdfstring{$g_1$}{g\_1}} \label{Sec:EstG1}

\begin{figure}[htbp]
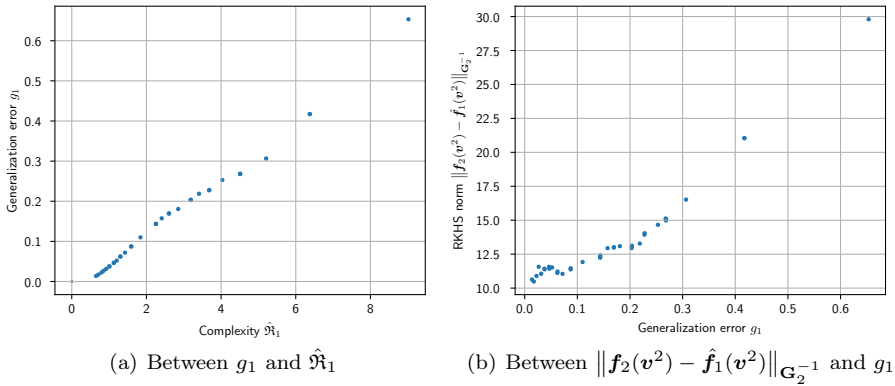

\subfigure[Between $g_1$ and $\hat{\mathfrak{R}}_1$]{\draftbox{\trimbox{4pt 2pt 12pt 16pt}{\scalebox{0.5}{\input{Figures/NLSE1DLine1.pgf}}}}}
\subfigure[Between $ \bnorm{\big}{ \bm{f}_2 \brbr{}{\bm{v}^2} - \hat{\bm{f}}_1 \brbr{}{\bm{v}^2} }_{\mathbf{G}_2^{-1}} $ and $g_1$]{\draftbox{\trimbox{-2pt 2pt 8pt 16pt}{\scalebox{0.5}{\input{Figures/NLSE1DLine2.pgf}}}}}
\caption{
Correlation between key quantities in the triangle {\color{blue!75!black}$\mathcal{G}_1$} and {\color{green!75!black}$\mathcal{C}_2$} of \cref{Fig:NoteEst}.
Detailed discussion on the computation is presented in \cref{Supp:Sec:Growth}.
}
\label{Fig:R1G1}
\end{figure}

As mentioned in \cref{Fig:NoteEst}, we apply \cref{Thm:Mar} in the triangle {\color{blue!75!black}$\mathcal{G}_1$} to estimate the generalization error $g_1$.
For training at level $ l = 1 $, the main term in the right hand side of \cref{Eq:Mar} is the complexity $ \hat{\mathfrak{R}}_1 = 2 R \bnorm{\big}{ \bm{f}_1 \brbr{}{\bm{v}^1} - \hat{\bm{f}}_0 \brbr{}{\bm{v}^1} }_{\mathbf{G}_1^{-1}} / \sqrt{M_1} $.
In the example of non-linear Schr\"odinger equation, we try with different $M_1$ and plot the correlation between complexity $\hat{\mathfrak{R}}_1$ and the generalization error $g_1$ in \cref{Fig:R1G1} (a).

With the clear linear correlation in \cref{Fig:R1G1} (a) as empirical justification, we consider the \emph{a priori} estimator to $g_1$:
\begin{equation} \label{Eq:G1HatG1}
\hat{g}_1 = \frac{ 2 R c_1 }{\sqrt{M_1}},
\quad c_1 = \bnorm{\big}{ \bm{f}_1 \brbr{}{\bm{v}^1} - \hat{\bm{f}}_0 \brbr{}{\bm{v}^1} }_{\mathbf{G}_1^{-1}}, \quad
R = \max_{ v \in \opsupp \mathcal{D} } \sqrt{\norm{ \Theta \rbr{ v, v } }_2}
\end{equation}
For justification of \cref{Eq:G1HatG1} from \cref{Thm:Mar}, we need to choose an appropriate cutoff $B_1$ such that (1) $ B_1 \gg g_1 $, so that we can equate the left hand side of \cref{Eq:Mar} with $g_1$ and (2) $ B_1 \ll c_1 $, so that the last two terms in \cref{Eq:Mar} are negligible.
We conclude from \cref{Fig:R1G1} (a) that in practice $g_1$ decreases when $M_1$ grows (actually in speed $ O \brbr{}{ 1 / \sqrt{M_1} } $).
However, $c_1$ does not decrease (but actually increases as pointed out in \cref{Supp:Sec:Growth}), so such $B_1$ satisfying $ g_1 \ll B_1 \ll c_1 $ exists for large $M_1$.
Actually, in the example of non-linear Schr\"odinger equation, $ c_1 = \num{56.543} $ when $ M_1 = 16 $ (averaged over 64 collections of samples $ \bcbr{}{ \brbr{}{ v^1_m, u^1_m } }_{ m = 1 }^{M_1} $) and is much greater than $g_1$ in \cref{Fig:R1G1} (a), so such $B_1$ exists even for a small $M_1$.
We note that the coefficients $c_1$ and $R$ can be practically computed without training the regression network.

\subsubsection{Estimating \texorpdfstring{$ \bnorm{\big}{ \bm{f}_l \brbr{}{\bm{v}^l} - \hat{\bm{f}}_{ l - 1 } \brbr{}{\bm{v}^l} }_{\mathbf{G}_l^{-1}} $}{|| f\_l (v\^{}l) - \^{}f\_(l-1) (v\^{}l) ||} for \texorpdfstring{$ 2 \le l \le L $}{2 <= l <= L}}

As mentioned in \cref{Fig:NoteEst}, in the triangle {\color{green!75!black}$\mathcal{C}_l$}, we need to combine the generalization error $ g_{ l - 1 } $ at the previous level and the difference between training samples $ u^l_m - f_{ l - 1 } \brbr{}{v^l_m} = f_l \brbr{}{v^l_m} - f_{ l - 1 } \brbr{}{v^l_m} $ to estimate the RKHS norm $ \bnorm{\big}{ \hat{f}_l - \hat{f}_{ l - 1 } }_{\mathcal{H}} = \bnorm{\big}{ \bm{f}_l \brbr{}{\bm{v}^l} - \hat{\bm{f}}_{ l - 1 } \brbr{}{\bm{v}^l} }_{\mathbf{G}_l^{-1}} $, which is critical to estimate $g_l$ as in the case of $g_1$.
In the example of non-linear Schr\"odinger equation, we plot the RKHS norm $ \bnorm{\big}{ \bm{f}_2 \brbr{}{\bm{v}^2} - \hat{\bm{f}}_1 \brbr{}{\bm{v}^2} }_{\mathbf{G}_2^{-1}} $ at level $ l = 2 $ with different $g_1$ in \cref{Fig:R1G1} (b) by varying the number of training samples $M_1$ at level $ l = 1 $.

We can notice a linear correlation with positive intercept in \cref{Fig:R1G1} (b).
Only when $ g_{ l - 1 } $ is large does $ \bnorm{\big}{ \bm{f}_l \brbr{}{\bm{v}^l} - \hat{\bm{f}}_{ l - 1 } \brbr{}{\bm{v}^l} }_{\mathbf{G}_l^{-1}} $ deviate from the intercept.
But more frequently $ g_{ l - 1 } $ is small because we are able to generate much more samples at coarse grids as mentioned in \cref{Sec:Back}.
In this case, we may assume $ \hat{f}_{ l - 1 } \approx f_{ l - 1 } $ and therefore $ \bnorm{\big}{ \bm{f}_l \brbr{}{\bm{v}^l} - \hat{\bm{f}}_{ l - 1 } \brbr{}{\bm{v}^l} }_{\mathbf{G}_l^{-1}} \approx \bnorm{\big}{ \bm{f}_l \brbr{}{\bm{v}^l} - \bm{f}_{ l - 1 } \brbr{}{\bm{v}^l} }_{\mathbf{G}_l^{-1}} \gg \bnorm{\big}{ \bm{f}_{ l - 1 } \brbr{}{\bm{v}^l} - \hat{\bm{f}}_{ l - 1 } \brbr{}{\bm{v}^l} }_{\mathbf{G}_l^{-1}} $.
To be more precise, we make the following first-order expansion:
\begin{multline} \label{Eq:Asy}
\bnorm{\big}{ \bm{f}_l \brbr{}{\bm{v}^l} - \hat{\bm{f}}_{ l - 1 } \brbr{}{\bm{v}^l} }_{\mathbf{G}_l^{-1}} = \bnorm{\big}{ \brbr{\big}{ \bm{f}_l \brbr{}{\bm{v}^l} - \bm{f}_{ l - 1 } \brbr{}{\bm{v}^l} } + \brbr{\big}{ \bm{f}_{ l - 1 } \brbr{}{\bm{v}^l} - \hat{\bm{f}}_{ l - 1 } \brbr{}{\bm{v}^l} } }_{\mathbf{G}_l^{-1}} \\
\approx \sqrt{\vphantom{\bpbr{\big}{}}\smash{ \bnorm{\big}{ \bm{f}_l \brbr{}{\bm{v}^l} - \bm{f}_{ l - 1 } \brbr{}{\bm{v}^l} }_{\mathbf{G}_l^{-1}}^2 + 2 \bpbr{\big}{ \bm{f}_{ l - 1 } \brbr{}{\bm{v}^l} - \hat{\bm{f}}_{ l - 1 } \brbr{}{\bm{v}^l}, \bm{f}_l \brbr{}{\bm{v}^l} - \bm{f}_{ l - 1 } \brbr{}{\bm{v}^l} }_{\mathbf{G}_l^{-1}} }} \\
\approx \bnorm{\big}{ \bm{f}_l \brbr{}{\bm{v}^l} - \bm{f}_{ l - 1 } \brbr{}{\bm{v}^l} }_{\mathbf{G}_l^{-1}} + \frac{ \bpbr{\big}{ \bm{f}_{ l - 1 } \brbr{}{\bm{v}^l} - \hat{\bm{f}}_{ l - 1 } \brbr{}{\bm{v}^l}, \bm{f}_l \brbr{}{\bm{v}^l} - \bm{f}_{ l - 1 } \brbr{}{\bm{v}^l} }_{\mathbf{G}_l^{-1}} }{\bnorm{\big}{ \bm{f}_l \brbr{}{\bm{v}^l} - \bm{f}_{ l - 1 } \brbr{}{\bm{v}^l} }_{\mathbf{G}_l^{-1}}},
\end{multline}
where
\begin{multline}
\label{Eq:Dual}
\bpbr{\big}{ \bm{f}_{ l - 1 } \brbr{}{\bm{v}^l} - \hat{\bm{f}}_{ l - 1 } \brbr{}{\bm{v}^l}, \bm{f}_l \brbr{}{\bm{v}^l} - \bm{f}_{ l - 1 } \brbr{}{\bm{v}^l} }_{\mathbf{G}_l^{-1}}
\\
\le \rbr{ \frac{1}{M_l} \bnorm{\big}{ \bm{f}_{ l - 1 } \brbr{}{\bm{v}^l} - \hat{\bm{f}}_{ l - 1 } \brbr{}{\bm{v}^l} }_{ 2, 1 } } \rbr{ M_l \bnorm{\big}{ \mathbf{G}_l^{-1} \brbr{\big}{ \bm{f}_l \rbr{\bm{v}^l} - \bm{f}_{ l - 1 } \brbr{}{\bm{v}^l} } }_{ 2, \infty } }.
\end{multline}
Here the $ \rbr{ 2, p } $-norm takes 2-norm on the spatial axis before taking $p$-norm over the index $m$ of training samples.
We note that
\begin{equation}
\frac{1}{M_l} \bnorm{\big}{ \bm{f}_{ l - 1 } \brbr{}{\bm{v}^l} - \hat{\bm{f}}_{ l - 1 } \brbr{}{\bm{v}^l} }_{ 2, 1 } = \frac{1}{M_l} \sum_{ m = 1 }^{M_l} \bnorm{\big}{ f_{ l - 1 } \brbr{}{v^l_m} - \hat{f}_{ l - 1 } \brbr{}{v^l_m} }_2
\end{equation}
is a Monte Carlo estimation of the generalization error $ g_{ l - 1 } $ due to the independence of samples $ \bcbr{}{ \brbr{}{ v^{ l - 1 }_m, u^{ l - 1 }_m } }{}_{ m = 1 }^{M_{ l - 1 }} $ at level $ l - 1 $ and $ \bcbr{}{ \brbr{}{ v^l_m, u^l_m } }{}_{ m = 1 }^{M_l} $ at level $l$.
In this way, we justify the estimation
\begin{gather}
\label{Eq:GLEst}
\bnorm{\big}{ \bm{f}_l \brbr{}{\bm{v}^l} - \hat{\bm{f}}_{ l - 1 } \brbr{}{\bm{v}^l} }_{\mathbf{G}_l^{-1}} \lesssim c_l + \frac{ d_l g_{ l - 1 } }{c_l}, \\
\label{Eq:GLCoef}
\quad c_l = \bnorm{\big}{ \bm{f}_l \brbr{}{\bm{v}^l} - \bm{f}_{ l - 1 } \brbr{}{\bm{v}^l} }_{\mathbf{G}_l^{-1}}, \quad d_l = M_l \bnorm{\big}{ \mathbf{G}_l^{-1} \brbr{\big}{ \bm{f}_l \brbr{}{\bm{v}^l} - \bm{f}_{ l - 1 } \brbr{}{\bm{v}^l} } }_{ 2, \infty }.
\end{gather}
We note that the coefficients $c_l$ and $d_l$ can be computed without training the network.
We further study their growth with respect to $M_l$ in \cref{Supp:Sec:Growth} empirically.

\subsubsection{Estimating \texorpdfstring{$g_l$}{g\_l} for \texorpdfstring{$ 2 \le l \le L $}{2 <= l <= L}}

\begin{figure}[htbp]
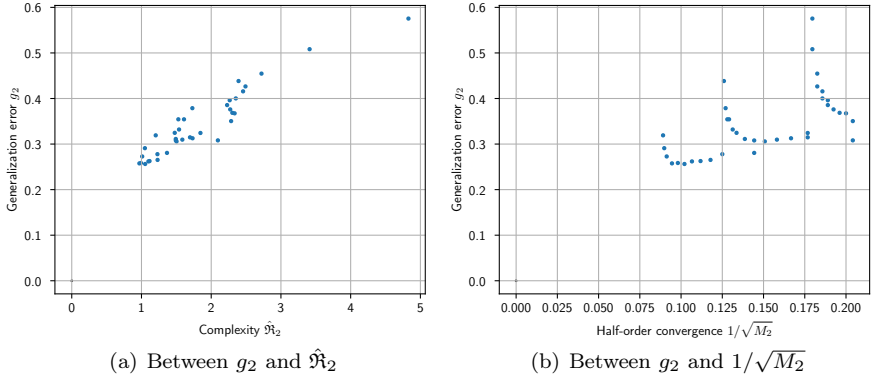

\subfigure[Between $g_2$ and $\hat{\mathfrak{R}}_2$]{\draftbox{\trimbox{4pt 2pt 12pt 16pt}{\scalebox{0.5}{\input{Figures/NLSE1DLine3.pgf}}}}}
\subfigure[Between $g_2$ and $ 1 / \sqrt{M_2} $]{\draftbox{\trimbox{4pt 2pt 12pt 16pt}{\scalebox{0.5}{\input{Figures/NLSE1DLine4.pgf}}}}}
\caption{
Correlation between the generalization error $g_2$ and the complexity $\hat{\mathfrak{R}}_2$ in the triangle {\color{blue!75!black}$\mathcal{G}_2$} of \cref{Fig:NoteEst}.
Detailed discussion on the computation is presented in \cref{Supp:Sec:NTKComp}. We note the importance of the RKHS norm $ \bnorm{\big}{ \bm{f}_2 \brbr{}{\bm{v}^2} - \hat{\bm{f}}_1 \brbr{}{\bm{v}^2} }_{\mathbf{G}_2^{-1}} $: if we leave it out in $\hat{\mathfrak{R}}_2$ and turn to consider $ 1 / \sqrt{M_2} $, we observe less correlation in (b) compared to (a).
}
\label{Fig:RLGL}
\end{figure}

As mentioned in \cref{Fig:NoteEst}, we apply \cref{Thm:Mar} in the triangle {\color{blue!75!black}$\mathcal{G}_l$} to estimate the generalization error $g_l$ from the previously discussed RKHS norm $ \bnorm{\big}{ \hat{f}_l - \hat{f}_{ l - 1 } }_{\mathcal{H}} = \bnorm{\big}{ \bm{f}_l \brbr{}{\bm{v}^l} - \hat{\bm{f}}_{ l - 1 } \brbr{}{\bm{v}^l} }_{\mathbf{G}_l^{-1}} $.
We display the correlation between the complexity $ \hat{\mathfrak{R}}_2 = 2 R \bnorm{\big}{ \bm{f}_2 \brbr{}{\bm{v}^2} - \hat{\bm{f}}_1 \brbr{}{\bm{v}^2} }_{\mathbf{G}_2^{-1}} / \sqrt{M_2}  $ and the generalization error $g_2$ at level $ l = 2 $ in the example of non-linear Schr\"odinger equation in \cref{Fig:RLGL}.

We can still witness the linear correlation in \cref{Fig:RLGL} (a) during fine-tuning (cf. \cref{Fig:R1G1} (a)).
Similar to the case of $g_1$, using the estimation to the RKHS norm \cref{Eq:GLEst}, we consider the \emph{a priori} estimator to $g_l$ with coefficients \cref{Eq:GLCoef}:
\begin{equation} \label{Eq:GL}
\hat{g}_l = \frac{ 2 R }{\sqrt{M_l}} \rbr{ c_l + \frac{ d_l \hat{g}_{ l - 1 } }{c_l} }.
\end{equation}

For the non-linear Schr\"odinger equation problem, we consider MLFT with $ L = 2 $ levels as mentioned in \cref{Sec:Back}.
The \emph{a priori} estimators to the generalization error in approximating the coarse-grid and fine-grid solvers at level $ l = 1 $ and $ l = 2 $ are respectively
\begin{equation} \label{Eq:G2Comp}
\hat{g}_1 = \frac{ 2 R c_1 }{\sqrt{M_1}}, \quad \hat{g}_2 = \frac{ 2 R }{\sqrt{M_2}} \rbr{ c_2 + \frac{ 2 d_2 R c_1 }{ c_2 \sqrt{M_1} } },
\end{equation}
where $c_1$, $R$ are computed from \cref{Eq:G1HatG1} and $c_2$, $d_2$ computed from \cref{Eq:GLCoef}.

\subsection{\emph{A posteriori} error estimator} \label{Sec:APost}

It is well known that \emph{a priori} error estimation may be pessimistic, especially in estimating the generalization of neural networks \cite{neyshabur_exploring_2017,arora_stronger_2018}.
Moreover, poor estimation may give poor solutions to the budget distribution problem \cref{Eq:Prog}.
Additionally, validity of the infinite-width NTK theory is questioned in finite-width scenarios which are mostly the case in practice \cite{allen-zhu_learning_2019,hanin_finite_2020}.
Hence, we may try to build practical \emph{a posteriori} error estimator to the generalization error $g_L$ in approximating the finest-grid solver $f_L$.
Combining theoretical insights and empirical observations, the \emph{a posteriori} error estimator may provide useful information about the generalization error for the budget distribution problem \cref{Eq:Prog}.
The \emph{a posteriori} error estimator may also extend to finite-width cases which are not solidly covered by the infinite-width NTK theory.
In detail, we substitute the coefficients $R$, $c_l$ and $d_l$ in \cref{Eq:G1HatG1,Eq:GL} by variables $a_l$ and $b_l$:
\begin{align} \label{Eq:GLPost}
\hat{g}_1 = \frac{a_1}{\sqrt{M_1}}, \quad \hat{g}_l = \frac{a_l}{\sqrt{M_l}} \rbr{ b_l + \hat{g}_{ l - 1 } }, \quad 2 \le l \le L.
\end{align}

We find values of the variables $a_l$ and $b_l$ \emph{a posteriori} from trials performing MLFT with different combinations of $M_l$.
By introducing a validation set, we compute the generalization error $g_l$ at the trials.
One particular method is the least square method: we fit $a_l$ and $b_l$ by minimizing the MSE between $ \log g_L $ and $ \log \hat{g}_L $ over trials.
The least square method needs at least $L$ trials to determine the function $ \hat{g}_L \rbr{ M_1, M_2, \cdots, M_l } $.
We also propose a heuristic method using only one trial of MLFT by setting $ 2 L - 1 $ equations: we set $ b_l = \ope_{ v \sim \mathcal{D} } \bnorm{\big}{ f_l \rbr{v} - f_{ l - 1 } \rbr{v} }_2 $ and $ g_l = \hat{g}_l $.
The motivation to set $b_l$ follows from the inequality in the triangle {\color{green!75!black}$\mathcal{C}_l$} of \cref{Fig:NoteEst} that $ \bnorm{\big}{ f_l \brbr{}{v} - \hat{f}_{ l - 1 } \brbr{}{v} }_2 \le \bnorm{\big}{ f_l \brbr{}{v} - f_{ l - 1 } \brbr{}{v} }_2 + \bnorm{\big}{ f_{ l - 1 } \brbr{}{v} - \hat{f}_{ l - 1 } \brbr{}{v} }_2 $.
For ML2MC \cite{lye_multi-level_2020}, the \emph{a posteriori} error estimator has the form $ \hat{g}_L = \sum_{ l = 1 }^L a_l / \sqrt{M_l} $. The variables $a_l$ can also be found by the heuristic method: we set $L$ equations equating the computed and the estimated generalization error at each level.

In the case of the non-linear Schr\"odinger equation problem, we have $ L = 2 $ and
\begin{equation}
\hat{g}_2^{\mathrm{MLFT}} \rbr{ M_1, M_2 } = \frac{a_2}{\sqrt{M_2}} \rbr{ b_2 + \frac{a_1}{\sqrt{M_1}} }, \quad \hat{g}_2^{\mathrm{ML2MC}} \rbr{ M_1, M_2 } = \frac{a_1}{\sqrt{M_1}} + \frac{a_2}{\sqrt{M_2}}.
\end{equation}
After determining values of the variables, one may solve the budget distribution problem \cref{Eq:Prog} to optimize the number of training samples $M_l$ at each level.
Here the optimal $M_l$ for ML2MC can be directly found by the Cauchy--Schwarz inequality, while the problem for MLFT is a little harder.
However, we note that $ 1 / \sqrt{\vphantom{\bpbr{\big}{}}\smash{ \prod_{ l' = 1 }^l M_{l'} }} $ is convex with respect to $ \rbr{ M_1, M_2, \cdots, M_L } $ for $ 1 \le l \le L $, so highly efficient solvers can be deployed to solve the problem \cref{Eq:Prog}.

\section{Numerical results} \label{Sec:Num}

We conduct numerical experiments to show the performance of MLFT and compare between algorithms.
We tackle problems with oscillations, discontinuities, or rough coefficients.
For these problems, the finest-grid solvers contain much detail.
Hence, it poses challenges for the regression network to have a small generalization error when there is a limited budget for generating training samples.

We first implement the neural networks in PyTorch \cite{paszke_pytorch_2019} and then migrate to JAX \cite{bradbury_jax_2018} and Stax \cite{bradbury_stax_2018} in order to use the package Neural Tangents \cite{novak_neural_2020} to calculate NTK.
The codes of experiments were run on Nvidia Tesla K80.
We separately generate three sets of samples for training, validation, and testing respectively according to the procedure described in \cref{Sec:Level}.
The generalization errors to construct the \emph{a posteriori} error estimator are computed on the validation set, while the reported errors \cref{Eq:GLDef,Eq:GLNorm} in figures are computed on the testing set.
The batch size of 32 is used for all the experiments.
We fix the seed of random batch and random initialization for better comparison.
We train the network for \num{10000} iterations at each level, when we observe that the training error is smaller than the testing error (e.g., in \cref{Fig:Comp,Fig:BurgErr,Fig:Ellip2}).

\subsection{Non-linear Schr\"odinger equation with oscillatory potentials} \label{Sec:NLSE}

We first consider the example of one-dimensional non-linear Schr\"odinger equation, or namely the Gross--Pitaevskiii equation.
This equation relates to Bose--Einstein condensation \cite{anglin_bose--einstein_2002,bagnato_bose--einstein_2015} and receives emerging research attention \cite{ilan_band-edge_2010}.
As already introduced in \cref{Sec:Back}, we consider the equation \cref{Eq:NLSE1,Eq:NLSE2} on domain $ \Omega = \sbr{ 0, 1 } $ with periodic boundary condition.
The solution map we are interested in is the potential-ground state map.

We set the dispersion coefficient to be $ \beta = 100 $.
The distribution $\mathcal{D}$ of potential $v$ is generated by
\begin{equation} \label{Eq:NLSEGen}
v \rbr{x} = \sum_{ k = 1 }^K
A_k
\rbr{ \alpha + \cos \rbr{ 2 \spi \omega_k x + \phi_k } }
\frac{1}{ \sqrt{ 2 \spi } \sigma_k }
\exp \rbr{ -\frac{ d \rbr{ x, c_k }^2}{ 2 \sigma_k^2 } },
\quad x \in \sbr{ 0, 1 }
\end{equation}
where $ d \rbr{ x, y } = x - y - \fbr{ x - y + 1 / 2 } \in \srbr{ -1 / 2, 1 / 2 } $.
We set $ K = 4 $, $ \alpha = 0.1 $, and sample $ A_k \sim \mathcal{U} \sbr{ -400, -200 } $, $\omega_k \sim \mathcal{U}\sbr{ 40, 80 } $, $ 1 / \sigma_k \sim \mathcal{U} \sbr{ 10, 20 } $, and $ \phi_k \sim \mathcal{U} \sbr{ 0, 2 \spi } $ independently.
The ground state $u$ is calculated by the gradient flow solver described in \cite{bao_computing_2004} with time step $ \tau = 1 $.

As mentioned in \cref{Sec:Back}, we consider $ L = 2 $ levels with grid steps $ h_1 = 1 / N_1 = 1 / 40 $ and $ h_2 = 1 / N_2 = 1 / 320 $ respectively.
From the visualization of potential-ground state pair in \cref{Fig:Comp}, we notice that large high-frequency components of the potential $v$ produce small oscillations in the ground state $u$.
Hence, we select $ R_{ 2 \to 1 } : \mathcal{V}_2 \to \mathcal{V}_1 $ to be Fourier restriction operator since the frequency of oscillations in the potential $v$ exceeds the Nyquist frequency of the coarse grid.
On the coarse grid only low-frequency components of the potential $v$ remain, and we select $ I_{ 1 \to 2 } : \mathcal{U}_1 \to \mathcal{U}_2 $ to be the cubic interpolation operator.
We observe that running times of coarse-grid and fine-grid solvers approximately satisfy $ t_2 = 64 t_1 $, as discussed in \cref{Supp:Sec:CompTime}.

We use the MNN-$\mathcal{H}$ network \cite{fan_multiscale_2019-1} as the regression network.
We use a network with 6 branches where each branch is a 5-layer 160-channel sub-network as specified in \cref{Supp:Sec:ModelSpec}.
We use Momentum \cite{qian_momentum_1999} as the optimizer with momentum coefficient $0.975$ and learning rate $10$.
The theory of NTK is proved to be applicable for the Momentum optimizer \cite{lee_wide_2019}.
The learning rate is seemingly huge but is applicable in the NTK parameterization \cite{jacot_neural_2018}.

\begin{figure}[htbp]
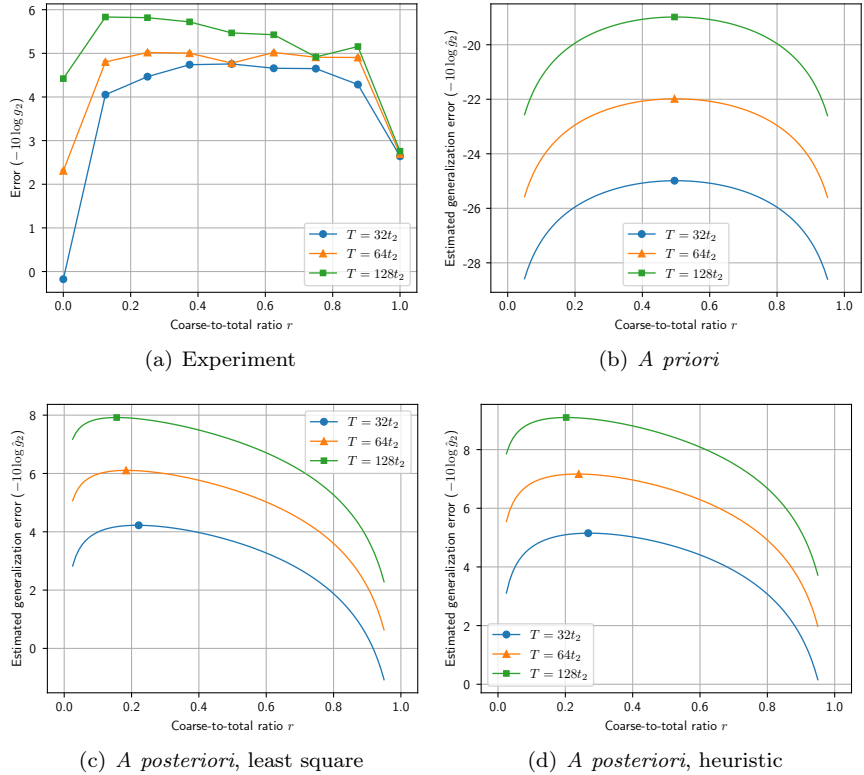

\centering
\subfigure[Experiment]{\draftbox{\trimbox{8pt 2pt 12pt 16pt}{\scalebox{0.5}{\input{Figures/NLSE1DWideExp.pgf}}}}}
\subfigure[\emph{A priori}]{\draftbox{\trimbox{4pt 2pt 12pt 16pt}{\scalebox{0.5}{\input{Figures/NLSE1DWidePri.pgf}}}}}
\subfigure[\emph{A posteriori}, least square]{\draftbox{\trimbox{8pt 2pt 12pt 16pt}{\scalebox{0.5}{\input{Figures/NLSE1DWidePostAll.pgf}}}}}
\subfigure[\emph{A posteriori}, heuristic]{\draftbox{\trimbox{8pt 2pt 12pt 16pt}{\scalebox{0.5}{\input{Figures/NLSE1DWidePostPart.pgf}}}}}
\caption
{
Error $g_2$ \cref{Eq:GLDef,Eq:GLNorm} and the estimated generalization error $\hat{g}_2$ for the non-linear Schr\"odinger equation problem.
In experiments (a) when $ r = 0 $ we use only fine-grid samples to train the regression network (i.e., single-level $ l = 2 $) and when $ r = 1 $ we use only coarse-grid samples (i.e., single-level $ l = 1 $).
We observe that MLFT with a selected $r$ is able to outperform single-level training with either $ l = 1 $ or $2$.
In particular, combination of $ T = 32 t_2 $ and $ r = 1 / 2 $ has smaller generalization error than $ T = 128 t_2 $ and $ r = 0 $.
This means that MLFT achieves the same generalization error while saving 75\% of the budget for generating training samples.
For error estimators, we observe that the \emph{a priori} error estimator (b) is pessimistic.
The \emph{a posteriori} error estimators capture the generalization error better.
The $r$ with the smallest estimated generalization error is marked, as the solution to the budget distribution problem \cref{Eq:Prog}.
A significant reduction of generalization error is still presented in experiments with these $r$ optimized from \emph{a posteriori} error estimators.
}
\label{Fig:NLSEErr}
\end{figure}

Recall that $M_1$ and $M_2$ are the numbers of training samples generated by the coarse-grid and fine-grid solvers $f_1$ and $f_2$ respectively.
We aim to reduce the generalization error $g_2$ under a fixed budget for generating training samples, i.e., $ M_1 t_1 + M_2 t_2 = T $.
To explore different combinations of $M_1$ and $M_2$ under a fixed budget $T$, we define the coarse-to-total ratio $ r = M_1 t_1 / \rbr{ M_1 t_1 + M_2 t_2 } $ to be the proportion of budget spent in generating coarse-grid samples.
We examine the relationship between the error $g_2$ \cref{Eq:GLDef,Eq:GLNorm} and the coarse-to-total ratio $r$ under a budget fixed $T$ in \cref{Fig:NLSEErr}.
We also plot the error estimators $\hat{g}_2$ to the generalization error introduced in \cref{Sec:APri,Sec:APost}.
The heuristic method to construct the \emph{a posteriori} error estimator runs a trial of MLFT with $ M_1 = 64 $ coarse-grid samples and $ M_2 = 32 $ fine-grid samples.

From the experiment, we observe that MLFT provides a significant reduction of the generalization error in approximating the fine-grid solver with a regression network.
The margin is larger when we are given smaller budget for generating training samples.
We notice that the $r$ with the smallest estimated generalization error decreases when the budget grows.
The \emph{a posteriori} error estimators correctly capture this tendency of decreasing $r$.

\subsection{Viscous Burgers' equation with discontinuous initial values} \label{Sec:Burg}

Burgers' equation plays an essential role in fluid dynamics and traffic flow modeling.
Since the viscous Burgers' equation contains both diffusion and non-linear convection terms, it is frequently taken as an example for conservation law solvers and model reduction methods \cite{lee_model_2020}.
The one-dimensional viscous Burgers' equation for $ u \rbr{ x; t } $ is
\begin{equation}
u_t + \rbr{\frac{u^2}{2}}_x = \kappa u_{ x x },
\quad x \in \Omega.
\end{equation}
We consider the equation on $ \Omega = \sbr{ 0, 1 } $ with a periodic boundary condition.
The solution map we are approximating is the time evolution map.
Formally, we are interested in the solution map $ v \mapsto u_v $ which takes initial value $ v = u \rbr{ \cdot; 0 } $ as parameter and gives terminal value $ u_v = u \rbr{ \cdot; t_{\mathrm{term}} } $ as solution.

We consider the problem with discontinuous initial values.
In detail, the distribution $\mathcal{D}$ of
parameter $v$ consists of step functions
\begin{equation}
v \rbr{x} = \sum_{ k = 1 }^K A_k \chi_{\srbr{ \rbr{ k - 1 } / K, k / K }} \rbr{x}, \quad x \in \sbr{ 0, 1 }.
\end{equation}
We set $K$ to be $40$ and sample $ A_k \sim \mathcal{N} \rbr{ 0, 1 } $ independently for $ 1 \le k \le K $.
We consider the solution at terminal time $ t_{\mathrm{term}} = 0.1 $ under diffusivity coefficient $ \kappa = 0.005 $.
For the numerical solver $F_l$ at each level, we use the first-order operator splitting scheme which includes Godnov scheme for convection and the explicit Euler scheme for diffusion.
Due to stability reasons, we choose the time step to be $ 10 h^2 $ on the grid with grid step $h$.
As a result, the time to run the level-$l$ solver satisfies $ t_l \propto 1 / h_l^3 $.
We choose the average operator to be $ R_{ L \to l } : \mathcal{V}_L \to \mathcal{V}_l $, and the linear interpolation operator to be $ I_{ l \to L } : \mathcal{U}_l \to \mathcal{U}_L $.

\begin{figure}[htbp]
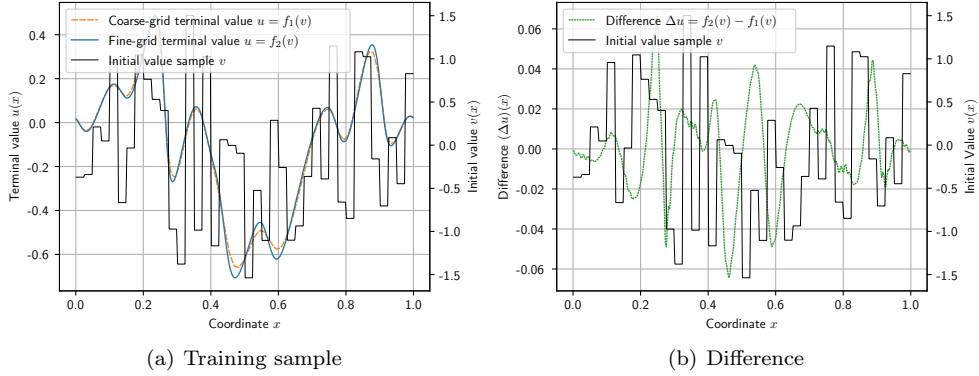

\centering
\subfigure[Training sample]{\draftbox{\trimbox{2pt 2pt -2pt 16pt}{\scalebox{0.5}{\input{Figures/Burg1DSGDDetFull.pgf}}}}}
\subfigure[Difference]{\draftbox{\trimbox{-2pt 2pt -2pt 16pt}{\scalebox{0.5}{\input{Figures/Burg1DSGDDetDiff.pgf}}}}}
\caption
{
Visualization of a training sample of the viscous Burgers' equation problem.
The coarse-grid sample suffers from numerical dissipation as shown in (a), which causes the difference between coarse-grid and fine-grid samples in (b).
}
\label{Fig:BurgDet}
\end{figure}

\begin{figure}[htbp]
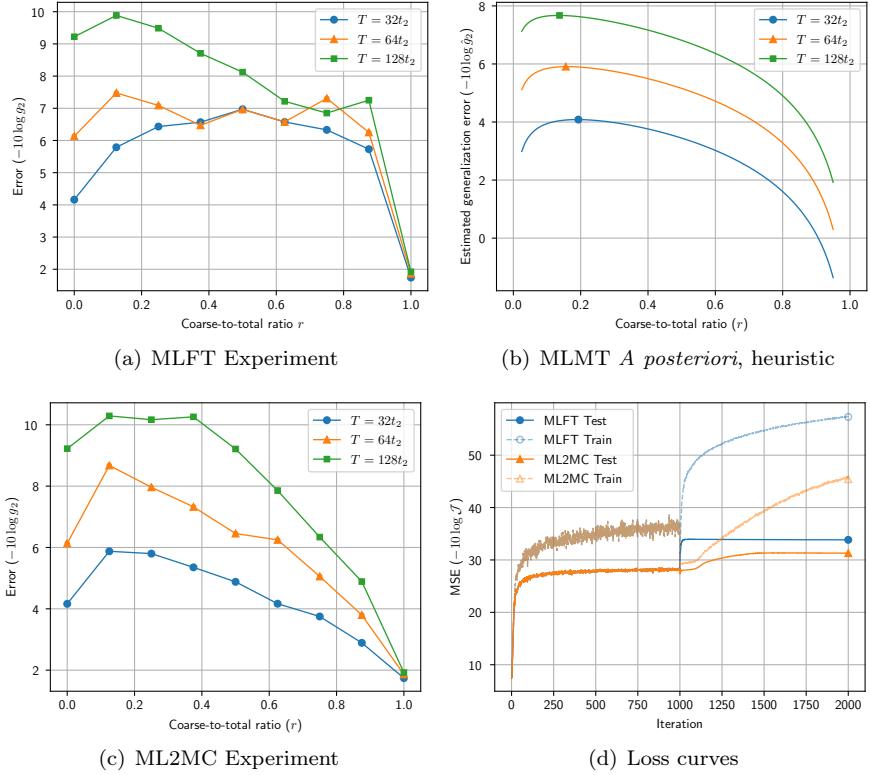

\centering
\subfigure[MLFT Experiment]{\draftbox{\trimbox{4pt 2pt 12pt 16pt}{\scalebox{0.5}{\input{Figures/Burg1DSGDExp.pgf}}}}}
\subfigure[MLMT \emph{A posteriori}, heuristic]{\draftbox{\trimbox{6pt 2pt 12pt 16pt}{\scalebox{0.5}{\input{Figures/Burg1DSGDPostPart.pgf}}}}}
\subfigure[ML2MC Experiment]{\draftbox{\trimbox{4pt 2pt 12pt 16pt}{\scalebox{0.5}{\input{Figures/Burg1DSGDExpKind2.pgf}}}}}
\subfigure[Loss curves]{\draftbox{\trimbox{4pt 2pt 12pt 16pt}{\scalebox{0.5}{\input{Figures/Burg1DSGDLoss.pgf}}}}}
\caption
{
Error $g_2$ \cref{Eq:GLDef,Eq:GLNorm} and the estimated generalization error $\hat{g}_2$ for the viscous Burgers' equation problem.
We observe that in experiments (a) the combination of $ T = 32 t_2 $ and $ r = 1 / 2 $ with MLFT has smaller generalization error than $ T = 64 t_2 $ and $ r = 0 $, achieving a 50\% saving of the budget.
The \emph{a posteriori} error estimator (b) captures the trend of error curves in (a).
In comparison between MLFT and ML2MC in (c), we observe that MLFT reduces the generalization error slightly better than ML2MC in the case $ T = 32 t_2 $, namely when the budget is small.
We plot the loss curves of MLFT and ML2MC in the case $ M_1 = 1024 $ and $ M_2 = 16 $ in (d), and find that at level $ l = 2 $ MLFT converges much faster (in less than 500 iterations) than ML2MC (in about 5000 iterations).
}
\label{Fig:BurgErr}
\end{figure}

We consider MLFT with $ L = 2 $ levels with grid steps $ h_1 = 1 / N_1 = 1 / 80 $ and $ h_2 = 1 / N_2 = 1 / 320 $.
In this case, $ t_2 = 64 t_1 $.
We visualize a training sample, namely an initial value-terminal value pair in \cref{Fig:BurgDet}.
We use a 9-layer convolutional neural network with 160 channels as the regression network, as specified in \cref{Supp:Sec:ModelSpec}.
We use the Momentum optimizer with learning rate $\num{10}$ and momentum coefficient $\num{0.975}$.
In this problem, MLFT still reduces the generalization error in approximating the fine-grid solver as shown in \cref{Fig:BurgErr}.
We compute the error $g_2$ \cref{Eq:GLDef,Eq:GLNorm} of MLFT and ML2MC with different coarse-to-total ratio $r$ under a fixed budget $T$.
We also plot the estimated generalization error $\hat{g}_2$ of the \emph{a posteriori} estimator constructed by the heuristic method which runs a trial of MLFT with $ M_1 = 64 $ coarse-grid samples and $ M_2 = 32 $ fine-grid samples.

From the experiments, we find that both MLFT and ML2MC reduces the generalization error $g_2$ in approximating the finest-grid solver $f_2$.
In the case that the budget for generating training samples is small, MLFT reduces the generalization error slightly better than ML2MC.
Besides, we find that at the fine level $ l = 2 $ ML2MC converges much slower than MLFT.
This is because of the effect of curriculum learning \cite{weinshall_curriculum_2018} as mentioned in \cref{Alg:TrainML2MC}.
In MLFT, the convergence of fine-tuning steps is boosted by the previous training or fine-tuning steps.
However, in ML2MC, networks are trained separately and converge separately.

\subsection{Diagonal of inverse to elliptic operators with rough coefficients}

The inverse to an elliptic operator is known as the Green function, which serves as the fundamental solution to elliptic PDE.
As an example, we are particularly interested in the diagonal of inverse $ u \rbr{\bm{x}} = g \rbr{ \bm{x}, \bm{x} } $ to the Schr\"odinger operator, namely $ g = \rbr{ -\Delta + v }^{-1} $.
This problem has been studied in \cite{lin_fast_2009,fan_bcr-net_2019} and has applications in quantum chemistry.
We aim to approximate the potential-diagonal of inverse map of the Schr\"odinger operator, which takes potential $ v \rbr{\bm{x}} $ as parameter and gives the diagonal of inverse $ u_v \rbr{\bm{x}} = u \rbr{\bm{x}} $ to the corresponding Schr\"odinger operator.

\begin{figure}[htbp]
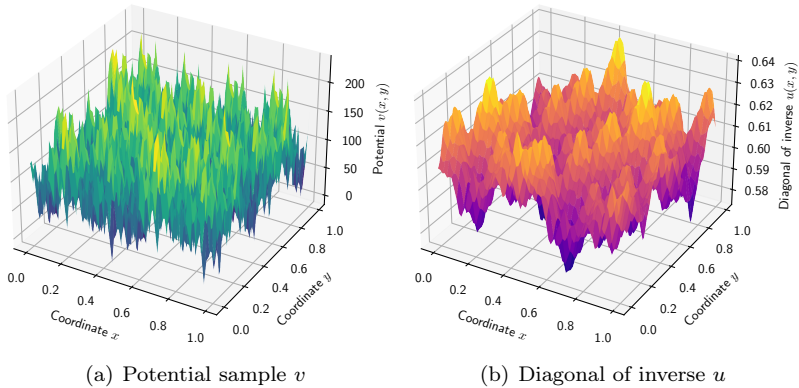

\centering
\subfigure[Potential sample $v$]{\draftbox{\trimbox{28pt 18pt 2pt 28pt}{\scalebox{0.5}{\input{Figures/Ellip2DRoughDetV.pgf}}}}}
\subfigure[Diagonal of inverse $u$]{\draftbox{\trimbox{28pt 18pt 2pt 28pt}{\scalebox{0.5}{\input{Figures/Ellip2DRoughDetG.pgf}}}}}
\caption
{
Visualization of a training sample of the diagonal of inverse problem.
We observe high-frequency oscillations in the potential $v$ (a), but the diagonal of the inverse $u$ to the corresponding Schr\"odinger operator (b) is smoother.
}
\label{Fig:EllipDet}
\end{figure}

Let us tackle the two-dimensional problem with $ \bm{x} = \rbr{ x, y } $ on $ \Omega = \sbr{ 0, 1 }^2 $ with a periodic boundary condition.
We generate rough potentials $v$ by
\begin{equation} \label{Eq:EllipVeq}
v \rbr{ x, y } = \sum_{ k = 1 }^K A_k \cos \rbr{ 2^k \spi \rbr{ x + y } + \phi_k } + \sum_{ k = 1 }^K B_k \cos \rbr{ 2^k \spi \rbr{ 2 x - y } + \psi_k } + C.
\end{equation}
For the distribution $\mathcal{D}$ of rough potentials, we set $ K = 6 $, $ C = 100 $ and sample $ A_k, B_k \sim \mathcal{N} \brbr{}{ 0, 20^2 } $ and $ \phi_k, \psi_k \sim \mathcal{U} \sbr{ 0, 2 \spi } $ independently.
We visualize a potential-diagonal of inverse pair in \cref{Fig:EllipDet}.
We invoke the SelInv algorithm \cite{lin_fast_2009} for the numerical solvers $F_l$, which takes $ O \brbr{}{ 1 / h^3 } $ times on a two-dimensional grid with grid step $ h = 1 / N $.
We choose $ R_{ L \rightarrow l } : \mathcal{V}_L \to \mathcal{V}_l $ to be the Fourier restriction operator, and $ I_{ l \rightarrow L } : \mathcal{U}_l \to \mathcal{U}_L $ to be the cubic interpolation operator.

\begin{figure}[htbp]
\centering
\draftbox{\trimbox{4pt 2pt 12pt 16pt}{\scalebox{0.5}{\input{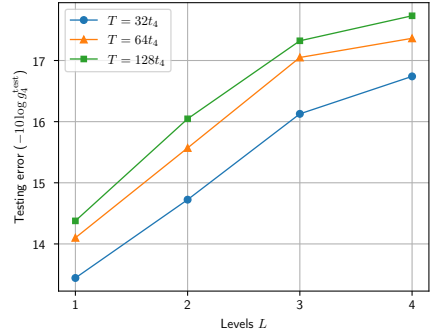}}}}
\caption{
Testing error $g_4^{\mathrm{test}}$ for the diagonal of inverse problem.
We observe that the error is significantly reduced under a fixed budget $T$ when the number of levels increases.
}
\label{Fig:Ellip1}
\end{figure}

We consider four different grids with grid step $ h_1 = 1 / N_1 = 1 / 10 $, $ h_2 = 1 / N_2 = 1 / 20 $, $ h_3 = 1 / N_3 = 1 / 40 $, and $ h_4 = 1 / N_4 = 1 / 80 $ respectively.
In this case, $ t_4 = 8 t_3 = 64 t_2 = 512 t_1 $.
We aim to approximate the finest-grid solver $f_4$ on the finest grid with grid step $ h_4 = 1 / 80 $.
We use a two-dimensional MNN-$\mathcal{H}$ network with 4 branches and 5 layers as the regression network, as specified in \cref{Supp:Sec:ModelSpec}.
We use the Adam optimizer with learning rate $\num{3.0e-4}$ \cite{kingma_adam_2015}.
Although the theory of NTK is not applicable in this case, we still consider the ability of MLFT in reducing the generalization error.
We perform MLFT with different number of levels $L$ for $ 1 \le L \le 4 $: we successively use training samples generated at the grid with grid step $ h_{ 4 + l - L } $ for $ l = 1, 2, \cdots, L $.
The number of training samples $M_l$ at different levels are optimized according to the budget distribution problem \cref{Eq:Prog}.
Here we use the \emph{a posteriori} error estimator constructed by the heuristic method which runs trials with $ M_1 = 256 $, $ M_2 = 128 $, $ M_3 = 64 $ and $ M_4 = 32 $.
We compute and report the testing error $g^{\mathrm{test}}_4$ with different numbers of levels $L$ under a fixed budget $T$ in \cref{Fig:Ellip1}.
We observe that as the number of levels grows, MLFT achieves smaller error in approximating the finest-grid solver under a fixed budget for generating training samples.

\begin{figure}[htbp]
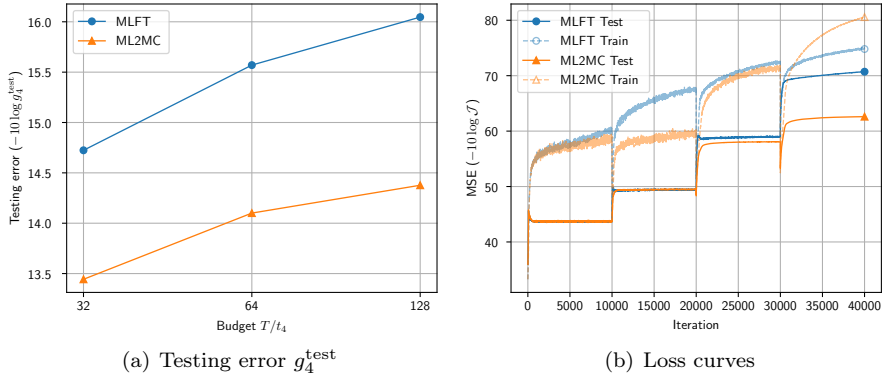

\centering
\subfigure[Testing error $g_4^{\mathrm{test}}$]{\draftbox{\trimbox{0pt 2pt 12pt 16pt}{\scalebox{0.5}{\input{Figures/EllipComp.pgf}}}}}
\subfigure[Loss curves]{\draftbox{\trimbox{4pt 2pt 12pt 16pt}{\scalebox{0.5}{\input{Figures/Ellip2DLoss.pgf}}}}}
\caption{
Comparison between MLFT and ML2MC for the diagonal of inverse problem.
In terms of accuracy (a), we observe that MLFT reduces the generalization error better than ML2MC.
In terms of convergence (b), we observe that MLFT converges faster than ML2MC at the level $ l = 3 $.
We note that the difference between MLFT and ML2MC becomes larger after entering level $ l = 4 $.
This is because MLFT corrects the generalization error in fine-tuning steps while ML2MC accumulates the generalization error.
The loss curves are extracted from the experiments in (a) under the fixed budget $ T = 32 t_4 $.
}
\label{Fig:Ellip2}
\end{figure}

We compare MLFT and ML2MC with $ L = 4 $ levels under different budgets $T$ for generating training samples in \cref{Fig:Ellip2}.
We compute the testing error $g^{\mathrm{test}}_4$ to compare the accuracy of MLFT and ML2MC.
We also plot the loss curves in order to understand the convergence and compare the efficiency.
For reference, we show the optimized number of training samples at different levels in \cref{Tbl:EllipWork} as the solution to the budget distribution problem \cref{Eq:Prog}.

\begin{table}[htbp]
\centering
\caption{
Optimized number of training samples at different levels for the diagonal of inverse problem.
The number of training samples is found by solving the budget distribution problem \cref{Eq:Prog} and is subject to $ M_1 t_1 + M_2 t_2 + M_3 t_3 + M_4 t_4 = T $ or $ M_1 / 512 + M_2 / 64 + M_3 / 8 + M_4 = T / t_4 $.
We observe that compared to MLFT, ML2MC generates much more coarse-grid samples in exchange of some fine-grid samples.
}
\label{Tbl:EllipWork}
\begin{tabular}{crrrrr}
\toprule
& $ T / t_4 $ & $M_1$ & $M_2$ & $M_3$ & $M_4$ \\
\midrule
\multirow{3}*{\makecell{MLFT \cref{Alg:TrainMLFT}
}}
& \num{32} & \num{53} & \num{51} & \num{46} & \num{25} \\
\cmidrule{2-6}
& \num{64} & \num{68} & \num{72} & \num{76} & \num{53} \\
\cmidrule{2-6}
& \num{128} & \num{86} & \num{101} & \num{126} & \num{111} \\
\midrule
\multirow{3}*{ML2MC \cref{Alg:TrainML2MC}} 
& \num{32} & \num{2143} & \num{653} & \num{54}  & \num{7} \\
\cmidrule{2-6}
& \num{64} & \num{4286} & \num{1305} & \num{107} & \num{13} \\
\cmidrule{2-6}
& \num{128} & \num{8572} & \num{2611} & \num{214} & \num{27} \\
\bottomrule
\end{tabular}
\end{table}

\section{Conclusion}

We establish the Multi-Level Fine-Tuning (MLFT) algorithm \cref{Alg:TrainMLFT} in this paper.
We aim to reduce the generalization error of regression networks in approximating solution maps and numerical solvers without spending more time in generating training samples.
Inspired by the resemblance between fine-tuning techniques and multi-scale methods, we train the regression networks on samples generated at the coarse grid and fine-tune on samples generated at finer grids.
In experiments, we achieve significant reduction of the generalization error with MLFT under a fixed budget for generating training samples, compared with base-lines trained with only coarse-grid or only fine-grid samples.
Thanks to the Neural Tangent Kernel (NTK) theory, we are able to perform analysis on the MLFT algorithm under assumptions.
Although it is notoriously difficult to find sharp bounds on the generalization error of neural networks, we make use of the form of generalization error provided by statistical machine learning theory and find the coefficients \emph{a posteriori}.
By minimizing the estimated generalization error given by the practical \emph{a posteriori} error estimator, we distribute the budget for generating training samples over levels.
In this way, we hope to provide practical guidance to reduce the generalization error with theoretical insight.

We have not exhaustively tested MLFT on many other intriguing numerical problems in scientific machine learning.
Problems with heavier scaling of time with respect to the grid step may enjoy greater reduction of generalization error with the help of MLFT.
It is also possible to generalize MLFT to other settings besides regression networks.
For the error estimators, we are interested in other methods to construct the \emph{a posteriori} error estimators besides the two mentioned in the paper.
More detailed analysis can be done, possibly providing new forms of the generalization error.

\section{Acknowledgements}

We would like to thank Tianle Cai for providing reference.
We also appreciate the discussion with Professor Bin Dong on the delivery of the subject.
We would like to express our gratitude to the anonymous reviewers together with Yiqing Xie, Weijie Chen, Jikai Hou, and Zeyu Zhao for careful proofreading and kind comments on the manuscript, which greatly improves the quality of the paper.

\bibliographystyle{siamplain}
\bibliography{Reference}

\makeatletter\@input{Supplement.aux.tex}\makeatother

\end{document}


\maketitle

\section{Proofs of theorems} \label{Sec:Pf}

\begin{proof}[Proof of \cref{Art:Thm:RadeKer}]
This theorem is a generalization of the Rademacher complexity of kernel classes \cite{bartlett_rademacher_2001,mohri_foundations_2018} to the vector-valued case \cite{maurer_rademacher_2006,sindhwani_scalable_2013}.

Firstly, we observe from the independence condition $ \nu_m \perp \sigma_m $ that
\begin{equation}
\begin{split}
\hat{\mathfrak{R}}_{ v, \nu } \brbr{\big}{ \mathcal{B} \brbr{\big}{ \hat{f}_{ l - 1 }, D } } &= \frac{1}{M_l} \ope_{ \sigma, \nu } \sup_{ \norm{f}_{\mathcal{H}} \le D } \sum_{ m = 1 }^{M_l} \sigma_m \nu_m^{\mathsf{T}} \brbr{\big}{ \hat{f}_{ l - 1 } + f } \brbr{}{v^l_m} \\
&= \frac{1}{M_l} \ope_{ \sigma, \nu } \sup_{ \norm{f}_{\mathcal{H}} \le D } \sum_{ m = 1 }^{M_l} \sigma_m \nu_m^{\mathsf{T}} f \brbr{}{v^l_m} \\
&= \hat{\mathfrak{R}}_{ v, \nu } \rbr{ \mathcal{B} \rbr{ 0, D } }.
\end{split}
\end{equation}
This means that the empirical Rademacher complexity of the $\mathcal{H}$-ball $ \mathcal{B} \brbr{\big}{ \hat{f}_{ l - 1 }, D } $ is not dependent on the center $ \hat{f}_{ l - 1 } $.

Moreover, we deduce
\begin{equation}
\begin{split}
\hat{\mathfrak{R}}_{ v^l, \nu } \brbr{\big}{ \mathcal{B} \brbr{\big}{ \hat{f}_{ l - 1 }, D } } &= \frac{1}{M_l} \ope_{ \sigma, \nu } \sup_{ \norm{f}_{\mathcal{H}} \le D } \sum_{ m = 1 }^{M_l} \bpbr{\big}{ f, \sigma_m \Theta \brbr{}{ \cdot, v^l_m } \nu_m }_{\mathcal{H}} \\
&= \frac{1}{M_l} \ope_{ \sigma, \nu } \sup_{ \norm{f}_{\mathcal{H}} \le D } \pbr{ f, \sum_{ m = 1 }^{M_l} \sigma_m \Theta \brbr{}{ \cdot, v^l_m } \nu_m }_{\mathcal{H}} \\
&= \frac{D}{M_l} \ope_{ \sigma, \nu } \norm{ \sum_{ m = 1 }^{M_l} \sigma_m \Theta \brbr{}{ \cdot, v^l_m } \nu_m }_{\mathcal{H}} \\
&\le \frac{D}{M_l} \sqrt{ \ope_{ \sigma, \nu } \norm{ \sum_{ m = 1 }^{M_l} \sigma_m \Theta \brbr{}{ \cdot, v^l_m } \nu_m }_{\mathcal{H}}^2 } \\
&= \frac{D}{M_l} \sqrt{ \sum_{ m = 1 }^{M_l} \ope_{\nu_m} \nu_m^{\mathsf{T}} \Theta \brbr{}{ v^l_m, v^l_m } \nu_m } \\
&\le \frac{ D R }{\sqrt{M_l}}
\end{split}
\end{equation}
as desired.
\end{proof}

\begin{proof}[Proof of \cref{Art:Thm:RadeLoss}]
Since we aim to approximate the level-$l$ solver function $f_l$ at level $l$, we define the loss function to be
\begin{equation} \label{Eq:LossFunc}
\ell \rbr{ u; v } = \norm{ f_l \rbr{v} - u }_2 \wedge B_l.
\end{equation}
We note the Lipchitz condition in $u$:
\begin{equation}
\ell \rbr{ u; v } - \ell \rbr{ u', v } \le \norm{ u - u' }_2.
\end{equation}

By leveraging the independence between $\sigma_1$ and $ \sigma_{ > 1 } = \rbr{ \sigma_2, \sigma_3, \cdots, \sigma_{M_l} } $, we have
\begin{equation}
\begin{split}
\hat{\mathfrak{R}}_{v^l} \brbr{\big}{ \mathcal{L} \brbr{\big}{ \hat{f}_{ l - 1 }, D; f_l, B_l } } 
&= \frac{1}{M_l} \ope_{\sigma} \sup_{ f \in \mathcal{B} \brbr{}{ \hat{f}_{ l - 1 }, D } } \sum_{ m = 1 }^{M_l} \sigma_m \ell \rbr{ f \brbr{}{v^l_m}; v^l_m } \\
&= \!
\begin{multlined}[t]
\frac{1}{ 2 M_l } \ope_{\sigma_{ > 1 }} \sup_{ f, f' \in \mathcal{B} \brbr{}{ \hat{f}_{l - 1 }, D } } \Bigg( \ell \brbr{\big}{ f \brbr{}{v^l_1}; v^l_1 } - \ell \brbr{\big}{ f' \brbr{}{v^l_1}; v^l_1 } \\
+ \sum_{ m = 2 }^{M_l} \sigma_m \ell \brbr{\big}{ f \brbr{}{v^l_m}; v^l_m } + \sum_{ m = 2 }^{M_l} \sigma_m \ell \brbr{\big}{ f' \brbr{}{v^l_m}; v^l_m } \Bigg).
\end{multlined}
\end{split}
\end{equation}
For a specific $ \sigma_{ > 1 } $, since
\begin{multline} \label{Eq:Supre}
\ell \brbr{\big}{ f \brbr{}{v^l_1}; v^l_1 } - \ell \brbr{\big}{ f' \brbr{}{v^l_1}; v^l_1 } 
\\
\le \bnorm{\big}{ f \brbr{}{v^l_1} - f' \brbr{}{v^l_1} }_2
= \sup_{ \norm{\nu_1}_2 = 1 } \nu_1^{\mathsf{T}} \brbr{\big}{ f \brbr{}{v^l_1} - f' \brbr{}{v^l_1} },
\end{multline}
we deduce
\begin{equation}
\begin{split}
\hat{\mathfrak{R}}_{v^l} \brbr{\big}{ \mathcal{L} \brbr{\big}{ \hat{f}_{ l - 1 }, D; f_l, B_l } } &\le \!
\begin{multlined}[t]
\frac{1}{ 2 M_l } \ope_{\sigma_{ > 1 }} \sup_{ f, f' \in \mathcal{B} \brbr{}{ \hat{f}_{ l - 1 }, D } } \sup_{ \norm{\nu_1}_2 = 1 } \Bigg( \nu_1^{\mathsf{T}} f \brbr{}{v^l_1} - \nu_1^{\mathsf{T}} f' \brbr{}{v^l_1} \\
+ \sum_{ m = 2 }^{M_l} \sigma_m \ell \brbr{\big}{ f \brbr{}{v^l_m}; v^l_m } + \sum_{ m = 2 }^{M_l} \sigma_m \ell \brbr{\big}{ f' \brbr{}{v^l_m}; v^l_m } \Bigg)
\end{multlined}
\\
&= \! 
\begin{multlined}[t]
\frac{1}{M_l} \sup_{\sarr{c}{ \norm{\nu_1}_2 \equiv 1 \\ \nu_1 \perp \sigma_1 }} \ope_{ \sigma, \nu_1 } \sup_{ f \in \mathcal{B} \brbr{}{ \hat{f}_{ l - 1 }, D } } \Bigg( \sigma_1 \nu_1^{\mathsf{T}} f \brbr{}{v^l_1} \\
+ \sum_{ m = 2 }^{M_l} \sigma_m \ell \brbr{\big}{ f \brbr{}{v^l_m}; v^l_m } \Bigg).
\end{multlined}
\end{split}
\end{equation}
Here the second supremum on $\nu_1$ is taken over all random vectors $\nu_1$ satisfying the normalization condition $ \norm{\nu_1}_2 \equiv 1 $ and the independence condition $ \nu_1 \perp \sigma_1 $.
This is because in the first supremum \cref{Eq:Supre} $\nu_1$ only depends on $ \sigma_{ > 1 } $ but not $\sigma_1$.

By repeating this process, we deduce
\begin{equation}
\begin{split}
\hat{\mathfrak{R}}_{v^l} \brbr{\big}{ \mathcal{L} \brbr{\big}{ \hat{f}_{ l - 1 }, D; f_l, B_l } } 
&\le \frac{1}{M_l} \sup_{\sarr{c}{ \norm{\nu_m} \equiv 1 \\ \nu_m \perp \sigma_m }} \ope_{ \sigma, \nu } \sup_{ f \in \mathcal{B} \brbr{}{ \hat{f}_{ l - 1 }, D } } \sum_{ m = 1 }^{M_l} \sigma_m \nu_m^{\mathsf{T}} f \brbr{}{v^l_m} \\
&\le \sup_{\sarr{c}{ \norm{\nu_m} \equiv 1 \\ \nu_m \perp \sigma_m }} \hat{\mathfrak{R}}_{ v^l, \nu } \brbr{\big}{ \mathcal{B} \brbr{\big}{ \hat{f}_{ l - 1 }, D } } \\
&\le \frac{ D R }{\sqrt{M_l}},
\end{split}
\end{equation}
where the last inequality follows from \cref{Art:Thm:RadeKer}.
\end{proof}

\begin{proof}[Proof of \cref{Art:Thm:Mar}]
Since $ \bnorm{\big}{ \hat{f}_l - \hat{f}_{ l - 1 } }_{\mathcal{H}} = \bnorm{\big}{ \bm{f}_l \brbr{}{\bm{v}^l} - \hat{\bm{f}}_{ l - 1 } \brbr{}{\bm{v}^l} }_{\mathbf{G}_l^{-1}} $ depends on the random training samples $ \bcbr{}{ \brbr{}{ v^l_m, u^l_m } }{}_{ m = 1 }^{M_l} $, we cannot directly assume a hypothesis space $ \mathcal{B} \brbr{\big}{ \hat{f}_{ l - 1 }, D } $ where $D$ is a deterministic constant.

We use the technique exploited in \cite{arora_fine-grained_2019}.
Since by \cref{Art:Thm:Proj}
\begin{equation} \label{Eq:Limit}
\bnorm{\big}{ \hat{f}_l - \hat{f}_{ l - 1 } }_{\mathcal{H}} = \bnorm{\big}{ \mathrm{P}_{\mathcal{M}_l} \brbr{\big}{ f_l - \hat{f}_{ l - 1 } } }_{\mathcal{H}} \le \bnorm{\big}{ f_l - \hat{f}_{ l - 1 } }_{\mathcal{H}},
\end{equation}
where the right hand side is known before random samples are generated, we may cover the possible hypothesis space by several $\mathcal{H}$-balls.

In detail, we assume $ \bgbr{\big}{ \bnorm{\big}{ f_l - \hat{f}_{ l - 1 } }_{\mathcal{H}} / B_l } = K_l $.
Hence, by applying \cref{Art:Thm:RadeMar} for $K_l$ times, with probability $ 1 - \delta $ on generating training samples $ \bcbr{}{ \brbr{}{ v^l_m, u^l_m } }{}_{ m = 1 }^{M_l} $, for all $ 1 \le k \le K_l $ it holds that
\begin{multline} \label{Eq:Union}
\sup_{\sarr{c}{ f \in \mathcal{B} \brbr{}{ \hat{f}_{ l - 1 }, k B_l } \\ f \brbr{}{v^l_m} = f_l \brbr{}{v^l_m} }} \ope_{ v^l \sim \mathcal{D} } \brbr{\big}{ \bnorm{\big}{ f_l \brbr{}{v^l} - f \brbr{}{v^l} }_2 \wedge B_l } \\
\le 2 \hat{\mathfrak{R}}_{v^l} \brbr{\big}{ \mathcal{L} \brbr{\big}{ \hat{f}_{ l - 1 }, k B_l; f_l, B_l } } + 3 B_l \sqrt{\frac{ \log K_l / \delta }{ 2 M_l }}.
\end{multline}
For a specific set of training samples $ \bcbr{}{ \brbr{}{ v^l_m, u^l_m } }{}_{ m = 1 }^{M_l} $, \cref{Eq:Limit} guarantees the existence of $k_0$ such that $ 1 \le k_0 \le K_l $ and
\begin{equation}
\rbr{ k_0 - 1 } B_l \le \bnorm{\big}{ \hat{f}_l - \hat{f}_{ l - 1 } }_{\mathcal{H}} \le k_0 B_l.
\end{equation}
As a result, $ \hat{f}_l \in \mathcal{B} \brbr{\big}{ \hat{f}_{ l - 1 }, k_0 B_l } $.
Since by \cref{Art:Thm:Proj} $ \hat{f}_l \brbr{}{v^l_m} = f_l \brbr{}{v^l_m} $, \cref{Eq:Union} can be applied to $ f = \hat{f}_l $ with $ k = k_0 $.
In this case,
\begin{equation} \label{Eq:RadeSp}
\hat{\mathfrak{R}}_{v^l} \brbr{\big}{ \mathcal{L} \brbr{\big}{ \hat{f}_{ l - 1 }, k_0 B_l; f_l, B_l } } \le \frac{ k_0 B_l R }{\sqrt{M_l}} \le \frac{ R \bnorm{\big}{ \hat{f}_l - \hat{f}_{ l - 1 } }_{\mathcal{H}} }{\sqrt{M_l}} + \frac{ B_l R }{\sqrt{M_l}}.
\end{equation}
Plugging \cref{Eq:RadeSp} into \cref{Eq:Union}, we obtain
\begin{equation}
\ope_{ v^l \sim \mathcal{D} } \brbr{\big}{ \bnorm{\big}{ f_l \brbr{}{v^l} - \hat{f}_l \brbr{}{v^l} }_2 \wedge B } \le \frac{ 2 R \bnorm{\big}{ \hat{f}_l - \hat{f}_{ l - 1 } }_{\mathcal{H}} }{\sqrt{M_l}} + \frac{ 2 B_l R }{\sqrt{M_l}} + 3 B_l \sqrt{\frac{ \log K_l / \delta }{ 2 M_l }}
\end{equation}
as desired.
\end{proof}

\section{Model specification} \label{Sec:ModelSpec}

\subsection{MNN-\texorpdfstring{$\mathcal{H}$}{H} network} \label{Sec:MNNH}

In the numerical experiments, we mainly use the MNN-$\mathcal{H}$ network \cite{fan_multiscale_2019-1}.
Since we are tackling problems with discontinuities, dilated convolution \cite{yu_multi-scale_2016} is not applicable because it may introduce greater jumps by sub-sampling with strides.
We make some adaptations to the MNN-$\mathcal{H}$ architecture for our numerical experiments.

We describe different branches of the
MNN-$\mathcal{H}$ network by the number of nodes.
For example, for a one-dimensional MNN-$\mathcal{H}$ network, the branch of 40 nodes has intermediate activations of shape $ \sbr{ \text{\texttt{batch\_size}}, 40 } $.
We leave out the adjacent component $A^{\rbr{\mathrm{ad}}}$ specified in \cite{fan_multiscale_2019-1} because it is covered by other branches.

We enlarge the restriction and interpolation maps, namely \textsf{LCR} and \textsf{LCI} layers in \cite{fan_multiscale_2019-1} to convolutional layers whose window sizes are larger than one.
We also initialize these layers with preset numerical restriction and interpolation operators.
This technique is found to stabilize the training and increase the accuracy.

We use the Neural Tangent Kernel (NTK) parameterization \cite{jacot_neural_2018}, which scales the output (activation) or each layer by $ 1 / \sqrt{C} $ where $C$ is the number of channels in the layer.
Except for manually initialized layers, parameters (weights and biases) of other layers are initialized independently using the unit Gaussian $ \mathcal{N} \rbr{ 0, 1 } $.
We always set the scaling coefficients for weights and biases to 1 ($\beta$ in \cite{jacot_neural_2018} and \texttt{W\_std}, \texttt{b\_std} in \cite{novak_neural_2020}).
We describe the parameterization and initialization of \textsf{LCR} and \textsf{LCI} layers.
We take the case of one dimension as an example.
We assume the spatial dimension of input to be $ N = N_L $, or equivalently the shape of input of the neural network to be $ \sbr{ \text{\texttt{batch\_size}}, N } $.
We consider the branch with $N_{\mathrm{sub}}$ nodes.
Weights of the \textsf{LCR} layer are initialized as
\begin{equation} \label{Eq:LCR}
W_{ c d w } = \sqrt{\frac{C}{ N / N_{\mathrm{sub}} }} W_c R_{ d w },
\end{equation}
where $ 1 \le c \le C $ is for output channels, $ 1 \le d \le N / N_{\mathrm{sub}} $ is for input channels, and $ -W \le w \le W $ is for the window of convolution.
Here $ R_{ d w } $ is a numerical restriction operator with normalization condition $ \sum_{ d = 1 }^{ N / N_{\mathrm{sub}} } \sum_{ w = -W }^W R_{ d w } = 1 $.
At initialization, we sample $W_c$ independently from $ \mathcal{N} \rbr{ 0, 1 } $ for $ 1 \le c \le C $.
Weights of the \textsf{LCI} layer are initialized as
\begin{equation} \label{Eq:LCI}
W_{ d c w } = \sqrt{\frac{C}{ N / N_{\mathrm{sub}} }} W_c I_{ d w },
\end{equation}
where $ 1 \le c \le C $ is for input channels, $ 1 \le d \le N / N_{\mathrm{sub}} $ is for output channels, and $ -W \le w \le W $ is for the window.
Here $ I_{ d w } $ is a numerical interpolation operator with normalization condition $ \sum_{ w = -W }^W I_{ d w } = 1 $.
At initialization, we sample $W_c$ independently from $ \mathcal{N} \rbr{ 0, 1 } $ for $ 1 \le c \le C $.
Biases of the \textsf{LCR} and \textsf{LCI} layers are sampled from $ \mathcal{N} \rbr{ 0, 1 } $ independently for initialization.

\begin{figure}[htbp]
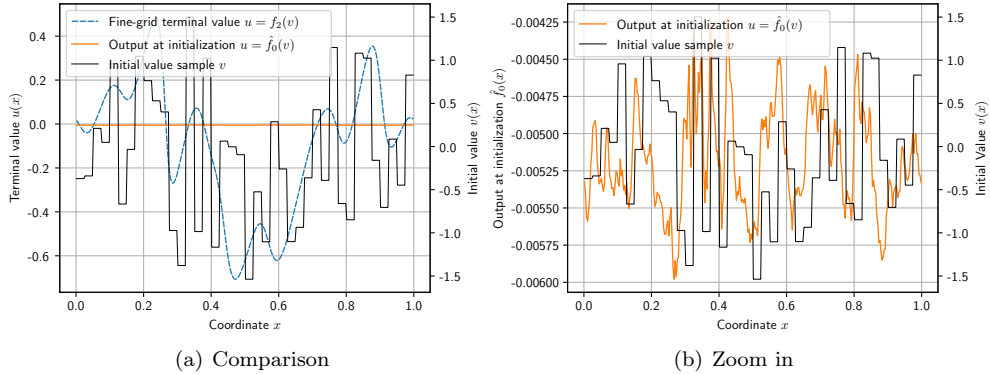

\centering
\subfigure[Comparison]{\draftbox{\trimbox{8pt 2pt -2pt 16pt}{\scalebox{0.5}{\input{Figures/Burg1DSGDInit1.pgf}}}}}
\subfigure[Zoom in]{\draftbox{\trimbox{-6pt 2pt 4pt 16pt}{\scalebox{0.5}{\input{Figures/Burg1DSGDInit2.pgf}}}}}
\caption
{
Visualization of the output of neural network at initialization for the viscous Burgers' equation problem.
We note that from (a) the output at initialization $ \hat{f}^0 \rbr{v} $ nearly vanishes compared to the solution given by the fine-grid solver $ f_2 \rbr{v} $.
However, from (b) it does not really vanish so the back-propagation gradients do not diminish.
}
\label{Fig:Init}
\end{figure}

In order to remove the randomness in computing $ c_1 = \bnorm{\big}{ \bm{f}_1 \brbr{}{\bm{v}^1} - \hat{\bm{f}}_0 \brbr{}{\bm{v}^1} }_{\mathbf{G}_1^{-1}} $ in \cref{Art:Sec:APri}, we scale the output of neural network by a small constant $\gamma$.
In this way, we have $ \hat{f}_0 \brbr{}{v^1_m} \approx 0 $ and we are able to compute $ c_1 \approx \bnorm{\big}{ \bm{f}_1 \brbr{}{\bm{v}^1} }_{\mathbf{G}_1^{-1}} $ regardless of random initialization of the network.
We note that $\gamma$ is not permitted to be too small, or otherwise the back-propagation gradients diminish.
We do not indent to bring the network into the lazy training regime by introduce this scaling as in \cite{chizat_lazy_2019}.
Actually, the scaling $ 1 / \sqrt{C} $ in NTK parameterization already results in lazy training.
For visualization, we plot the output of neural network at initialization in \cref{Fig:Init}.

\subsection{Models for experiments}

For the non-linear Schr\"odinger equation problem, we use a one-dimensional MNN-$\mathcal{H}$ network with spatial dimension of input $ N = 320 $.
We deploy 6 branches with the number of nodes $ N_{\mathrm{sub}} = 10, 20, 40, 80, 160, $ and $ 320 $ respectively.
The numbers of layers and channels of each branch are 5 and 160 respectively.
We initialize the \textsf{LCR} layer with average operator and the \textsf{LCI} layer with linear interpolation operator.
The window size is 3 for \textsf{LCR} and \textsf{LCI} layers and 7 for intermediate convolutional layers.
The scaling constant is $ \gamma = \num{3.0e-4} $.

For the viscous Burgers' equation problem, we cast a Convolution Neural Network (CNN) as a one-dimensional MNN-$\mathcal{H}$ network with spatial dimension of input $ N = 320 $.
As a CNN, we deploy only 1 branch with the number of nodes $ N_{\mathrm{sub}} = 320 $.
There are 9 layers and 160 channels in the network.
We initialize the \textsf{LCR} layer with average operator and the \textsf{LCI} layer with linear interpolation operator.
The window size is 3 for \textsf{LCR} and \textsf{LCI} layers and 15 for intermediate convolutional layers.
The scaling constant is $ \gamma = \num{3.0e-2} $.

For the diagonal of inverse problem, we use a two-dimensional MNN-$\mathcal{H}$ network with spatial dimension of input $ N \times N = 80 \times 80 $.
We deploy 4 branches with the number of nodes $ N_{\mathrm{sub}} \times N_{\mathrm{sub}} = 10 \times 10, 20 \times 20, 40 \times 40 $ and $ 80 \times 80 $ respectively.
The numbers of layers and channels of each branch are 3 and 40 respectively.
We initialize the \textsf{LCR} layer with average operator and the \textsf{LCI} layer with linear interpolation operator.
The window size is $ 3 \times 3 $ for \textsf{LCR} and \textsf{LCI} layers and $ 7 \times 7 $ for intermediate convolutional layers.
The scaling constant is $ \gamma = \num{1.0e-3} $.

\section{Problem specification}

\subsection{Computation time} \label{Sec:CompTime}

\begin{figure}[htbp]
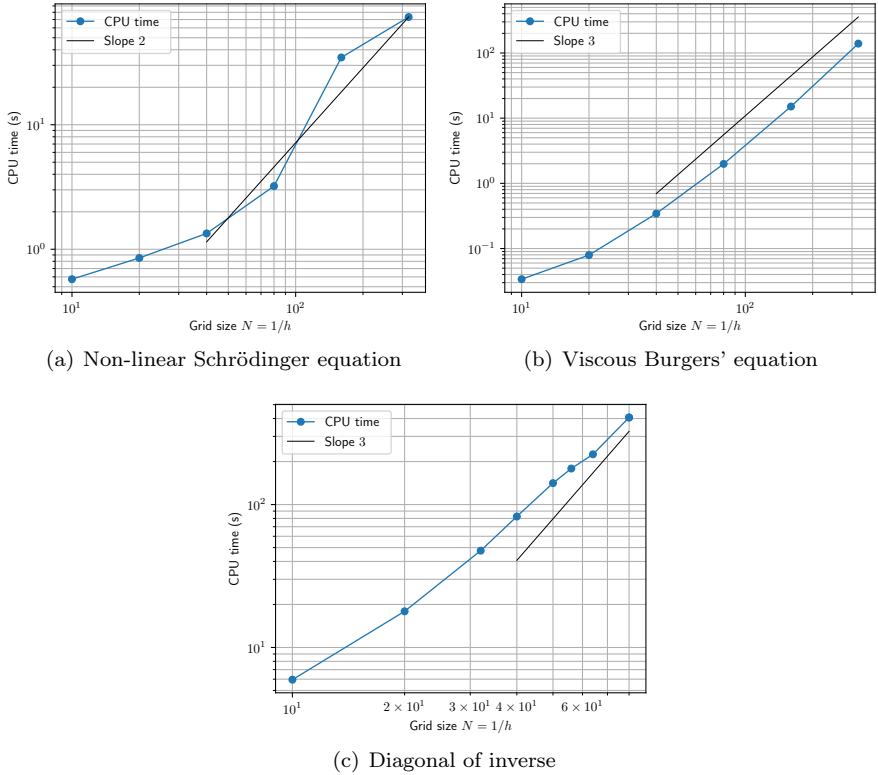

\centering
\subfigure[Non-linear Schr\"odinger equation]{\draftbox{\trimbox{4pt 2pt 12pt 16pt}{\scalebox{0.5}{\input{Figures/NLSE1DScale.pgf}}}}}
\subfigure[Viscous Burgers' equation]{\draftbox{\trimbox{2pt 2pt 12pt 16pt}{\scalebox{0.5}{\input{Figures/Burg1DScale.pgf}}}}}
\subfigure[Diagonal of inverse]{\draftbox{\trimbox{4pt 2pt 12pt 16pt}{\scalebox{0.5}{\input{Figures/Ellip1DScale.pgf}}}}}
\caption{
Measured CPU time in generating $ M = \num{16384} $ training samples at different grid steps.
For the non-linear Schr\"odinger equation problem, we observe a rapid increment of computation time for $ 80 \le N \le 160 $.
This is caused by the high-frequency oscillations \cref{Art:Eq:NLSEGen} which are only visible on fine grids and hinder the convergence of the gradient flow solver. 
For the viscous Burgers' equation problem we observe a scaling of order $ O \brbr{}{ 1 / h^3 } $.
For the diagonal of inverse problem, we also observe a scaling of order $ O \brbr{}{ 1 / h^3 } $ as given in \cite{lin_fast_2009}.
}
\label{Fig:Time}
\end{figure}

We measure the time $t_l$ to generate training samples at different levels.
For each problem in \cref{Art:Sec:Num}, we generate $ M = \num{16384} $ training samples with 6 cores in parallel on Intel Xeon CPU E5-2690.
We vary the grid step $ h = 1 / N $ and measure the CPU time for generating training samples.
The results are shown in \cref{Fig:Time}.
We note that the grid is of size $N$ for one-dimensional problems (the non-linear Schr\"odinger equation problem and the viscous Burgers' equation problem), and of size $ N \times N $ for two-dimensional problems (the diagonal of inverse problem).

According to the measurements, we are able to justify the following relationships between $t_l$, which are critical to the budget distribution problem \cref{Art:Eq:Prog}.
For the non-linear Schr\"odinger equation problem, we set $ t_2 = 64 t_1 $ as observed in \cref{Fig:Time}.
For the viscous Burgers' equation problem, we set $ t_2 = 64 t_1 $ according to the scaling law in time.
For the diagonal of inverse problem, we set $ t_4 = 8 t_3 = 64 t_2 = 512 t_1 $ according to the scaling law in time.

\subsection{Pre-processing}

We perform pre-processing on the generated training samples before we fit regression networks on the samples.
Pre-processing procedures are done for most machine learning algorithms (e.g., normalization, whitening, or centering).

For the non-linear Schr\"odinger equation problem, we observe that when the potential $ v \rbr{x} \equiv 0 $ the corresponding ground state is $ u \rbr{x} \equiv 1 $.
To center the ground states, we record $ u^l_m - 1 = f_l \brbr{}{v^l_m} - 1 $ instead of $u^l_m$ for the training sample pairs.

We do not perform pre-processing on the training samples for the viscous Burgers' equation problem.

For the diagonal of inverse problem, we notice that the potential \cref{Art:Eq:EllipVeq} has expectation $ \ope v^l \rbr{ x, y } \equiv C $.
To center the potentials, we record $ v^l_m - C $ instead of $v^l_m$ for the training samples pairs.
Since the two-dimensional Green function has singularities of logarithmic growth near the diagonal, further numerical treatment is required to match solutions generated at different grids.
We find from Fourier analysis that when the potential $ v^l \rbr{ x, y } \equiv C $, the diagonal of inverse to the discretized Schr\"odinger operator on the $ N_l \times N_l $ grid $\Omega_l$ is
\begin{equation}
u^l \rbr{ x, y } \equiv \tilde{C}_l = \sum_{ i = 1 }^{N_l} \sum_{ j = 1 }^{N_l} \frac{1}{ 4 N_l^2 \rbr{ \sin^2 i \spi / N_l + \sin^2 j \spi / N_l } + C }.
\end{equation}
To center the diagonal of inverses, we record $ u^l_m - \tilde{C}_l = f_l \brbr{}{v^l_m} - \tilde{C}_l $ for the training sample pairs when the grid step is $ h_l = 1 / N_l $.

\section{Neural Tangent Kernel} \label{Sec:NTKComp}

\subsection{Computation}

The computation of Conjugate Kernel (CK) and NTK of infinite-width networks has been extensively studied and efficiently implemented \cite{jacot_neural_2018,arora_exact_2019}.
We use the Neural Tangents package \cite{novak_neural_2020} together with JAX \cite{bradbury_jax_2018} and Stax \cite{bradbury_stax_2018} to compute the NTK.
Since the infinite-width NTK can be computed by sequentially applying CK and NTK operations corresponding to the layers \cite{novak_neural_2020}, we only need to specify the CK and NTK operations of manually initialized \textsf{LCR} and \textsf{LCI} layers introduce in \cref{Sec:MNNH}.

For \textsf{LCR} layers, the initialization \cref{Eq:LCR} can be recognized as restricting the input as a function on $ N \times N $ grid to $ N_{\mathrm{sub}} \times N_{\mathrm{sub}} $ grid with the restriction operator $R$ and then randomizing over channels.
As a result, the CK $\Sigma$ of a \textsf{LCR} layer is the inner product kernel of restricted functions, namely
\begin{equation}
\Sigma_{ n n' } \rbr{ v, v' } = \frac{1}{ \rbr{ 2 W + 1 } \rbr{ N / N_{\mathrm{sub}}} } \sum_{ w = -W }^W \sum_{ w' = -W }^W \sum_{ d = 1 }^{ N / N_{\mathrm{sub}} } \sum_{ d' = 1 }^{ N / N_{\mathrm{sub}} } R_{ d w } R_{ d' w' } v_{ n d } v'_{ n' d' } + 1,
\end{equation}
where $ 1 \le n, n' \le N_{\mathrm{sub}} $ are for spatial dimension and $ v, v' $ are recognized as shape $ N_{\mathrm{sub}} \times \rbr{ N / N_{\mathrm{sub}} } $ by reshaping the tensor.
The NTK of a manually initialized \textsf{LCR} layer is identical to a randomly initialized one, which turns out to be a standard convolutional layer.

For \textsf{LCI} layers, the CK $\Sigma$ is obtained by conjugating the CK $\Sigma'$ from the previous layer by the interpolation operator $I$ as
\begin{equation}
\Sigma_{ n n' d d' } \rbr{ v, v' } = \frac{1}{ 2 W + 1 } \sum_{ w = -W }^W \sum_{ w' = -W }^W I_{ d w } I_{ d' w' } \Sigma'_{ \rbr{ n + w } \rbr{ n' + w' } } \rbr{ v, v' } + \delta_{ d d' },
\end{equation}
where $ 1 \le n, n' \le N_{\mathrm{sub}} $ are for spatial dimension, $ 1 \le d, d' \le N / N_{\mathrm{sub}} $ are for output channels, and $ \delta_{ d d' } $ is the Kronecker delta.
The NTK $\Theta$ is also obtained by $\Sigma'$ and the NTK from the previous layer $\Theta'$
\begin{multline}
\Theta_{ n n' d d' } \rbr{ v, v' } = \frac{1}{ 2 W + 1 } \sum_{ w = -W }^W \sum_{ w' = -W }^W I_{ d w } I_{ d' w' } \Theta'_{ \rbr{ n + w } \rbr{ n' + w' } } \rbr{ v, v' } \\
+ \delta_{ d d' } + \frac{1}{ 2 W + 1 } \sum_{ w = -W }^W \Sigma'_{ \rbr{ n + w } \rbr{ n' + w } } \rbr{ v, v' } \delta_{ d d' },
\end{multline}
where $ 1 \le n, n' \le N_{\mathrm{sub}} $ are for spatial dimension and $ 1 \le d, d' \le N / N_{\mathrm{sub}} $ are for output channels.
The CK $\Sigma$ and NTK $\Theta$ above are of shape $ N_{\mathrm{sub}} \times N_{\mathrm{sub}} \times \rbr{ N / N_{\mathrm{sub}} } \times \rbr{ N / N_{\mathrm{sub}} } $, which can be recognized as shape $ \rbr{ N_{\mathrm{sub}} \times N / N_{\mathrm{sub}} } \times \rbr{ N_{\mathrm{sub}} \times N / N_{\mathrm{sub}} } = N \times N $ by reshaping the tensor.

\subsection{Convergence of Neural Tangent Kernel} \label{Sec:NTKConv}

\begin{figure}[htbp]
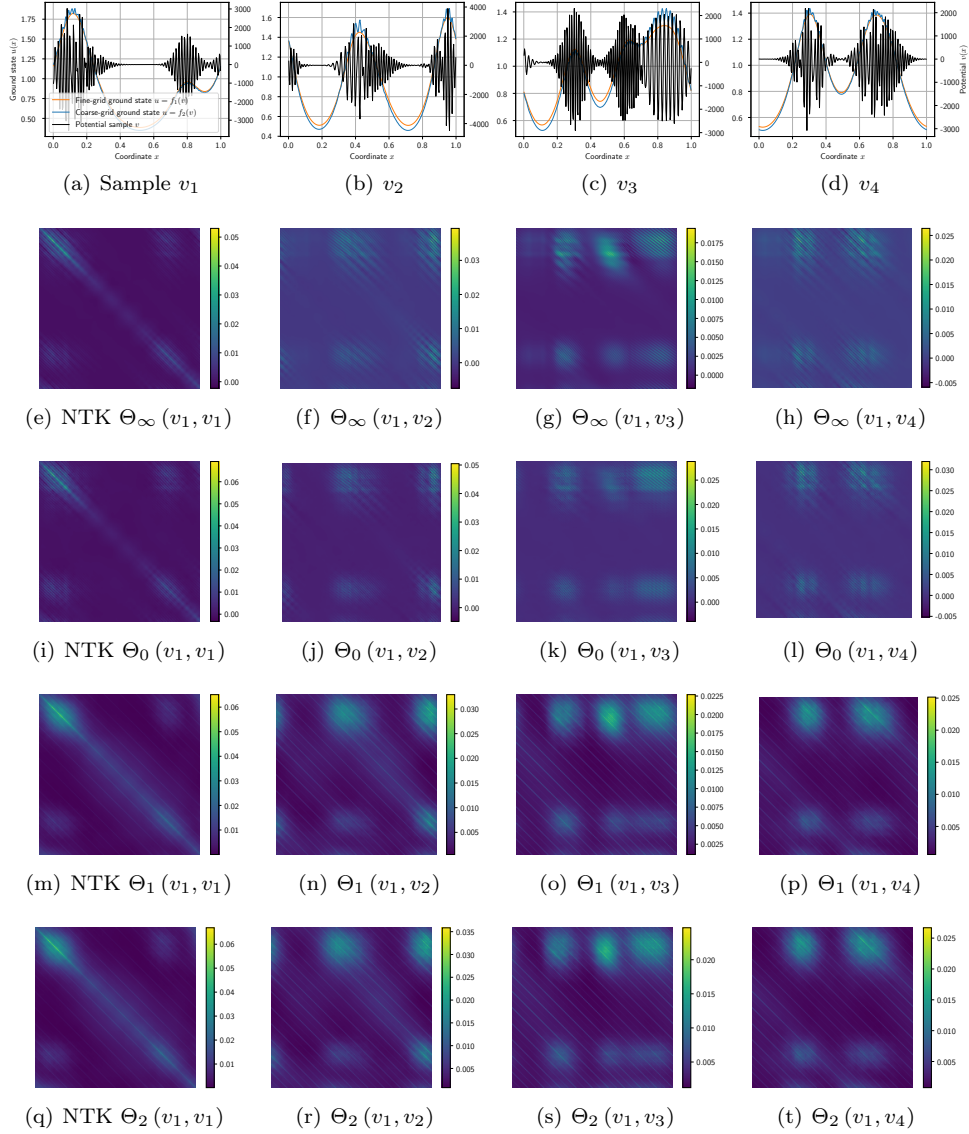

\centering
\subfigure[Sample $v_1$]{\draftbox{\trimbox{-2pt 0pt 0pt 6pt}{\scalebox{0.31}{\input{Figures/NLSE1DNTKPlot0.pgf}}}}}
\subfigure[$v_2$]{\draftbox{\trimbox{4pt 0pt 0pt 6pt}{\scalebox{0.31}{\input{Figures/NLSE1DNTKPlot1.pgf}}}}}
\subfigure[$v_3$]{\draftbox{\trimbox{4pt 0pt 0pt 6pt}{\scalebox{0.31}{\input{Figures/NLSE1DNTKPlot2.pgf}}}}}
\subfigure[$v_4$]{\draftbox{\trimbox{4pt 0pt -4pt 6pt}{\scalebox{0.31}{\input{Figures/NLSE1DNTKPlot3.pgf}}}}}
\subfigure[NTK $ \Theta_{\infty} \rbr{ v_1, v_1 } $]{\draftbox{\trimbox{-2pt 2pt 4pt 2pt}{\scalebox{0.31}{\input{Figures/NLSE1DNTKHeatI0.pgf}}}}}
\subfigure[$ \Theta_{\infty} \rbr{ v_1, v_2 } $]{\draftbox{\trimbox{-2pt 2pt 4pt 2pt}{\scalebox{0.31}{\input{Figures/NLSE1DNTKHeatI1.pgf}}}}}
\subfigure[$ \Theta_{\infty} \rbr{ v_1, v_3 } $]{\draftbox{\trimbox{-2pt 2pt 4pt 2pt}{\scalebox{0.31}{\input{Figures/NLSE1DNTKHeatI2.pgf}}}}}
\subfigure[$ \Theta_{\infty} \rbr{ v_1, v_4 } $]{\draftbox{\trimbox{-2pt 2pt 4pt 2pt}{\scalebox{0.31}{\input{Figures/NLSE1DNTKHeatI3.pgf}}}}}
\subfigure[NTK $ \Theta_0 \rbr{ v_1, v_1 } $]{\draftbox{\trimbox{-2pt 2pt 4pt 2pt}{\scalebox{0.31}{\input{Figures/NLSE1DNTKHeat00.pgf}}}}}
\subfigure[$ \Theta_0 \rbr{ v_1, v_2 } $]{\draftbox{\trimbox{-2pt 2pt 4pt 2pt}{\scalebox{0.31}{\input{Figures/NLSE1DNTKHeat01.pgf}}}}}
\subfigure[$ \Theta_0 \rbr{ v_1, v_3 } $]{\draftbox{\trimbox{-2pt 2pt 4pt 2pt}{\scalebox{0.31}{\input{Figures/NLSE1DNTKHeat02.pgf}}}}}
\subfigure[$ \Theta_0 \rbr{ v_1, v_4 } $]{\draftbox{\trimbox{-2pt 2pt 4pt 2pt}{\scalebox{0.31}{\input{Figures/NLSE1DNTKHeat03.pgf}}}}}
\subfigure[NTK $ \Theta_1 \rbr{ v_1, v_1 } $]{\draftbox{\trimbox{-2pt 2pt 4pt 2pt}{\scalebox{0.31}{\input{Figures/NLSE1DNTKHeat10.pgf}}}}}
\subfigure[$ \Theta_1 \rbr{ v_1, v_2 } $]{\draftbox{\trimbox{-2pt 2pt 4pt 2pt}{\scalebox{0.31}{\input{Figures/NLSE1DNTKHeat11.pgf}}}}}
\subfigure[$ \Theta_1 \rbr{ v_1, v_3 } $]{\draftbox{\trimbox{-2pt 2pt 4pt 2pt}{\scalebox{0.31}{\input{Figures/NLSE1DNTKHeat12.pgf}}}}}
\subfigure[$ \Theta_1 \rbr{ v_1, v_4 } $]{\draftbox{\trimbox{-2pt 2pt 4pt 2pt}{\scalebox{0.31}{\input{Figures/NLSE1DNTKHeat13.pgf}}}}}
\subfigure[NTK $ \Theta_2 \rbr{ v_1, v_1 } $]{\draftbox{\trimbox{-2pt 2pt 4pt 2pt}{\scalebox{0.31}{\input{Figures/NLSE1DNTKHeat40.pgf}}}}}
\subfigure[$ \Theta_2 \rbr{ v_1, v_2 } $]{\draftbox{\trimbox{-2pt 2pt 4pt 2pt}{\scalebox{0.31}{\input{Figures/NLSE1DNTKHeat41.pgf}}}}}
\subfigure[$ \Theta_2 \rbr{ v_1, v_3 } $]{\draftbox{\trimbox{-2pt 2pt 4pt 2pt}{\scalebox{0.31}{\input{Figures/NLSE1DNTKHeat42.pgf}}}}}
\subfigure[$ \Theta_2 \rbr{ v_1, v_4 } $]{\draftbox{\trimbox{-2pt 2pt 4pt 2pt}{\scalebox{0.31}{\input{Figures/NLSE1DNTKHeat43.pgf}}}}}
\caption
{
Visualization of NTK for the non-linear Schr\"odinger equation problem.
We plot the potential samples $v_m$ in (a-d), show infinite-width NTK $ \Theta_{\infty} \rbr{ v_1, v_m } $ in (e-h), and show finite-width NTK at initialization $ \Theta_0 \rbr{ v_1, v_m } $ in (i-l).
We also show finite-width NTK after one level of training $ \Theta_1 \rbr{ v_1, v_m } $ in (m-p) and after two levels of training $ \Theta_2 \rbr{ v_1, v_m } $ in (q-t).
We observe that finite-width NTK at initialization $\Theta_0$ is a qualitatively good approximation to infinite-width NTK.
Although finite-width NTK changes during the first level of training from $\Theta_0$ to $\Theta_1$, it still qualitatively captures prominent features of the input.
We observe that finite-width NTK does not change much during the second level of fine-tuning from $\Theta_1$ to $\Theta_2$.
}
\label{Fig:NTKPlot}
\end{figure}

We study the convergence of NTK in the example of non-linear Schr\"odinger equation.
We use the Neural Tangents package \cite{novak_neural_2020} to compute the infinite-width NTK $\Theta_{\infty}$.
Meanwhile, we can also use the package to compute the NTK of finite-width networks according to \cref{Art:Eq:NTKDef} directly.
We consider the NTK $\Theta_0$ of finite-width network at initialization, namely the NTK of $\hat{f}_0$.
We also consider finite-width NTK $\Theta_1$ of $\hat{f}_1$ after the first level of training with $ M_1 = 1024 $ coarse-grid samples and $\Theta_2$ of $\hat{f}_2$ after two levels of training, fine-tuned with $ M_2 = 16 $ fine-grid samples starting from $\hat{f}_1$.
We sample four $ v_m \sim \mathcal{D} $ for $ 1 \le m \le 4 $, and show the corresponding NTK $ \Theta_{\infty} \rbr{ v_1, v_m } $, $ \Theta_0 \rbr{ v_1, v_m } $, $ \Theta_1 \rbr{ v_1, v_m } $, and $ \Theta_2 \rbr{ v_1, v_m } $ in \cref{Fig:NTKPlot}.
We note that since the output of the network is multi-dimensional, the NTK is a matrix-valued kernel.

\begin{figure}[htbp]
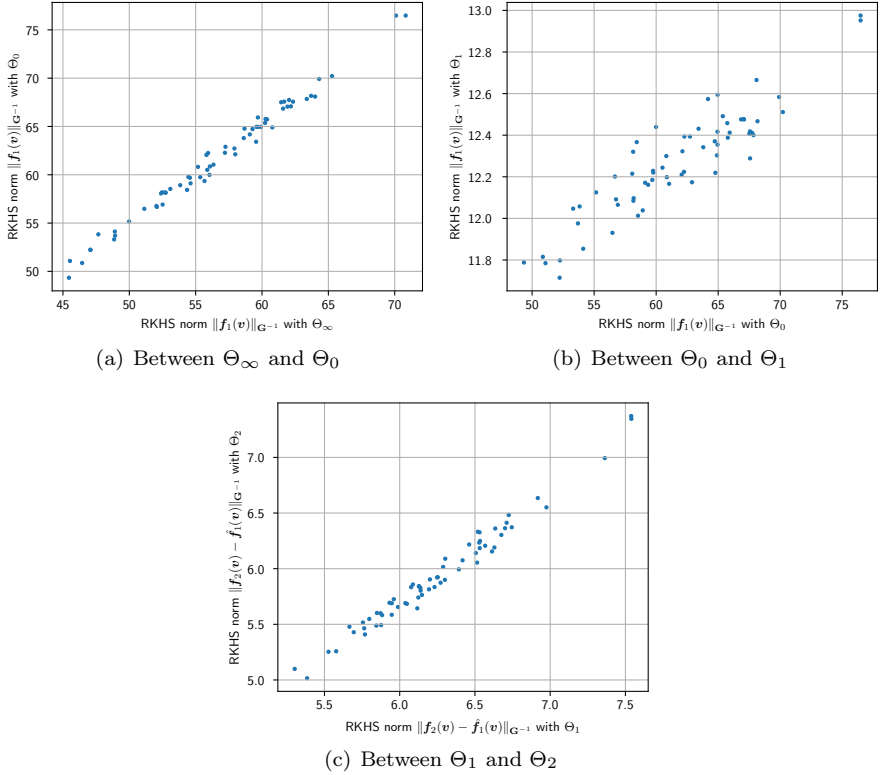

\centering
\subfigure[Between $\Theta_{\infty}$ and $\Theta_0$]{\draftbox{\trimbox{2pt 2pt 12pt 16pt}{\scalebox{0.5}{\input{Figures/NLSE1DFI0.pgf}}}}}
\subfigure[Between $\Theta_0$ and $\Theta_1$]{\draftbox{\trimbox{0pt 2pt 12pt 16pt}{\scalebox{0.5}{\input{Figures/NLSE1DF01.pgf}}}}}
\subfigure[Between $\Theta_1$ and $\Theta_2$]{\draftbox{\trimbox{0pt 2pt 12pt 16pt}{\scalebox{0.5}{\input{Figures/NLSE1DD14New.pgf}}}}}
\caption
{
Correlation between RKHS norms $ \bnorm{\big}{ \bm{f}_1 \brbr{}{\bm{v}} }_{\mathbf{G}^{-1}} $ computed with infinite-width NTK $\Theta_{\infty}$ and finite-width NTK $\Theta_0$, $\Theta_1$, and correlation between RKHS norms $ \bnorm{\big}{ \bm{f}_2 \brbr{}{\bm{v}} - \hat{\bm{f}}_1 \brbr{}{\bm{v}} }_{\mathbf{G}^{-1}} $ computed with finite-width NTK $\Theta_1$ and $\Theta_2$.
We observe a strong correlation between the RKHS norms with $\Theta_{\infty}$ and $\Theta_0$, which supports the convergence of NTK.
Although the RKHS norm with $\Theta_1$ changes in magnitude compared to $\Theta_0$, we still observe a linear correlation, which supports the construction of \emph{a posteriori} error estimator.
We again observe a strong correlation between the RKHS norms with $\Theta_1$ and $\Theta_2$, which supports the validity of NTK analysis in the fine-tuning steps.
}
\label{Fig:NTKCorr}
\end{figure}

Since the training process turns out to be kernel ``ridgeless'' regression in the NTK regime \cite{jacot_neural_2018}, we are particular interested in the RKHS norm $ \bnorm{\big}{ \hat{f}_1 - \hat{f}_0 }_{\mathcal{H}} = \bnorm{\big}{ \bm{f}_1 \brbr{}{\bm{v}^1} - \hat{\bm{f}}_0 \brbr{}{\bm{v}^1} }_{\mathbf{G}_1^{-1}} \approx \bnorm{\big}{ \bm{f}_1 \brbr{}{\bm{v}^1} }_{\mathbf{G}_1^{-1}} $ and $ \bnorm{\big}{ \hat{f}_2 - \hat{f}_1 }_{\mathcal{H}} = \bnorm{\big}{ \bm{f}_2 \brbr{}{\bm{v}^2} - \hat{\bm{f}}_1 \brbr{}{\bm{v}^2} }_{\mathbf{G}_2^{-1}} $ in $\mathcal{H}$ induced by the NTK.
For $ M = 16 $ samples $ \bcbr{}{v_m}{}_{ m = 1 }^{M} $, we compute the RKHS norm $ \bnorm{\big}{ \bm{f}_1 \brbr{}{\bm{v}}  }_{\mathbf{G}^{-1}} $ with infinite-width NTK $\Theta_{\infty}$ and finite-width NTK $\Theta_0$, $\Theta_1$ for the Gram matrix $\mathbf{G}$ of $ \bcbr{}{v_m}{}_{ m = 1 }^{M} $.
We also compute the RKHS norm $ \bnorm{\big}{ \bm{f}_2 \brbr{}{\bm{v}} - \hat{\bm{f}}_1 \brbr{}{\bm{v}} }_{\mathbf{G}^{-1}} $ with finite-width NTK $\Theta_1$ and $\Theta_2$.
We draw 64 collections of samples $ \bcbr{}{v_m}{}_{ m = 1 }^{M} $ by $ v_m \sim \mathcal{D} $ independently and study the correlation between RKHS norms with $\Theta_{\infty}$, $\Theta_0$, $\Theta_1$, and $\Theta_2$ in \cref{Fig:NTKCorr}.

\subsection{Growth of coefficients \texorpdfstring{$c_l$}{c\_l} and \texorpdfstring{$d_l$}{d\_l}} \label{Sec:Growth}

We introduce the coefficients $c_l$ for $ 1 \le l \le L $ and $d_l$ for $ 2 \le l \le L $ to construct the \emph{a priori} error estimator in \cref{Art:Sec:APri}.
According to \cref{Art:Eq:G1HatG1,Art:Eq:GLCoef}, computing $c_l$ and $d_l$ involves computing the Gram matrix $\mathbf{G}_l$ and its inverse.
However, we note that the Gram matrix is of size $ M_l N_L^d \times M_l N_L^d $, which becomes prohibitively huge when $M_l$ becomes large.
For example, we can only compute the Gram matrix with $ M_l = 24 $ samples at a time on Nvidia Tesla K80 for the non-linear Schr\"odinger equation problem where $ N_L = 320 $ and $ d = 1 $.
As a result, we only use $16$ samples instead of $M_1$ samples when computing the RKHS norms in \cref{Art:Fig:R1G1,Art:Fig:RLGL}.
In this way, the linear correlation shown in \cref{Art:Fig:R1G1} shows the convergence of the generalization error $g_1$ in speed $ O \rbr{ 1 / \sqrt{M_1} } $.

\begin{figure}[htbp]
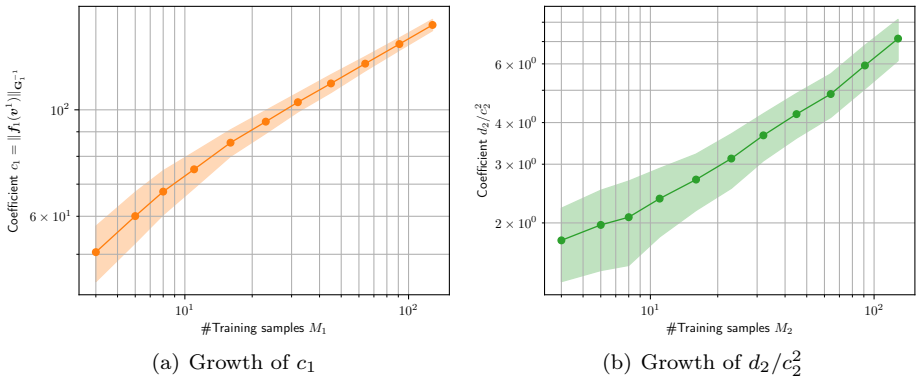

\centering
\subfigure[Growth of $c_1$]{\draftbox{\trimbox{-6pt 2pt 12pt 16pt}{\scalebox{0.5}{\input{Figures/NLSE1DGrowthC.pgf}}}}}
\subfigure[Growth of $ d_2 / c_2^2 $]{\draftbox{\trimbox{-4pt 2pt 12pt 16pt}{\scalebox{0.5}{\input{Figures/NLSE1DGrowthDCC.pgf}}}}}
\caption
{
Growth of $c_1$ and $ d_2 / c_2^2 $ with respect to $M_1$ and $M_2$ respectively.
We observe that the growth of $c_1$ and $ d_2 / c_2^2 $ are relatively slow.
The slope of $c_1$ with respect to $M_1$ in (a) is $ \num{0.30588} < 1 / 2 $, and the slope of $ d_2 / c_2^2 $ with respect to $M_2$ in (b) is $ \num{0.40594} < 1 $.
As a result, the growth exponents do not qualitatively change our result.
We generate 64 collections of samples, and the shaded area corresponds to the variance.
}
\label{Fig:NLSESmallGrow}
\end{figure}

We study the growth of the coefficients $c_l$ and $d_l$ with respect to the number of training samples $M_l$.
We still take the non-linear Schr\"odinger equation as an example.
For computing the Gram matrix with more samples, we consider a smaller problem: we use the coarse grid with grid step $ h_1 = 1 / N_1 = 1 / 20 $ and the fine grid with grid step $ h_2 = 1 / N_2 = 1 / 40 $ respectively.
We still use \cref{Art:Eq:NLSEGen} for the distribution $\mathcal{D}$ of potential $v$ but set $ K = 4 $ and sample $ A_k \sim \mathcal{U} \sbr{ -400, -200 } $, $ \omega_k \sim \mathcal{U} \sbr{ 8, 16 } $, $ 1 / \sigma_k \sim \mathcal{U} \sbr{ 4, 8 } $, and $ \phi_k \sim \mathcal{U} \sbr{ 0, 2 \spi } $ independently.
Besides the growth of $ c_1 = \bnorm{\big}{ \bm{f}_1 \brbr{}{\bm{v}^1} }_{\mathbf{G}_1^{-1}} $ with respect to $M_1$, since by \cref{Art:Eq:GLCoef}
\begin{equation}
\hat{g}_2 = \frac{ 2 R c_2 }{\sqrt{M_2}} \rbr{ 1 + \frac{ d_2 \hat{g}_1 }{c_2^2} },
\end{equation}
we study the growth of $ d_2 / c_2^2 = M_2 \bnorm{\big}{ \mathbf{G}_2^{-1} \brbr{\big}{ \bm{f}_2 \brbr{}{\bm{v}^2} - \bm{f}_1 \brbr{}{\bm{v}^2} } }_{ 2, \infty } / \bnorm{\big}{ \bm{f}_2 \brbr{}{\bm{v}^2} - \bm{f}_1 \brbr{}{\bm{v}^2} }_{\mathbf{G}_2^{-1}}^2 $ with respect to $M_2$.
The growth curves are presented in \cref{Fig:NLSESmallGrow}.
We observe that the growth of $c_1$ and $ d_2 / c_2^2 $ with respect to $M_1$ and $M_2$ respectively approximately follows the power law.
Due to the presence of growth exponents, it is possible to introduce exponent corrections for the \emph{a posteriori} error estimator.

\bibliographystyle{siamplain}
\bibliography{Reference}

\makeatletter\@input{Article.aux.tex}\makeatother